\documentclass[12pt]{article}%
\usepackage{indentfirst}
\usepackage[a4paper,left=3cm,right=2cm]{geometry}
\usepackage{amsfonts}
\usepackage{amssymb}
\usepackage{amsmath}
\usepackage{fancyhdr}
\usepackage{hyperref}
\usepackage{graphicx}%
\providecommand{\U}[1]{\protect\rule{.1in}{.1in}}
\providecommand{\U}[1]{\protect\rule{.1in}{.1in}}
\newtheorem{theorem}{Theorem}

\newtheorem{axiom}[theorem]{Axiom}

\newtheorem{conjecture}[theorem]{Conjecture}
\newtheorem{corollary}[theorem]{Corollary}

\newtheorem{definition}[theorem]{Definition}
\newtheorem{example}[theorem]{Example}
\newtheorem{exercise}[theorem]{Exercise}
\newtheorem{lemma}[theorem]{Lemma}

\newtheorem{proposition}[theorem]{Proposition}
\newtheorem{remark}[theorem]{Remark}

\newenvironment{proof}[1][Proof]{\noindent\textbf{#1.} }{\ \rule{0.5em}{0.5em}}

\makeatletter

\newcount\@hour\newcount\@minute\chardef\@x10\chardef\@xv60
\def\tcitime{
\def\@time{%
  \@minute\time\@hour\@minute\divide\@hour\@xv
  \ifnum\@hour<\@x 0\fi\the\@hour:%
  \multiply\@hour\@xv\advance\@minute-\@hour
  \ifnum\@minute<\@x 0\fi\the\@minute
  }}%

\def\x@hyperref#1#2#3{%
   \catcode`\~ = 12
   \catcode`\$ = 12
   \catcode`\_ = 12
   \catcode`\# = 12
   \catcode`\& = 12
   \y@hyperref{#1}{#2}{#3}%
}

\def\y@hyperref#1#2#3#4{%
   #2\ref{#4}#3
   \catcode`\~ = 13
   \catcode`\$ = 3
   \catcode`\_ = 8
   \catcode`\# = 6
   \catcode`\& = 4
}

\@ifundefined{hyperref}{\let\hyperref\x@hyperref}{}
\@ifundefined{msihyperref}{\let\msihyperref\x@hyperref}{}

\@ifundefined{qExtProgCall}{\def\qExtProgCall#1#2#3#4#5#6{\relax}}{}
\def\QCTOpt[#1]#2{%
  \def\QCTOptB{#1}
  \def\QCTOptA{#2}
}
\def\QCTNOpt#1{%
  \def\QCTOptA{#1}
  \let\QCTOptB\empty
}
\def\Qct{%
  \@ifnextchar[{%
    \QCTOpt}{\QCTNOpt}
}
\def\QCBOpt[#1]#2{%
  \def\QCBOptB{#1}%
  \def\QCBOptA{#2}%
}
\def\QCBNOpt#1{%
  \def\QCBOptA{#1}%
  \let\QCBOptB\empty
}
\def\Qcb{%
  \@ifnextchar[{%
    \QCBOpt}{\QCBNOpt}%
}
\def\PrepCapArgs{%
  \ifx\QCBOptA\empty
    \ifx\QCTOptA\empty
      {}%
    \else
      \ifx\QCTOptB\empty
        {\QCTOptA}%
      \else
        [\QCTOptB]{\QCTOptA}%
      \fi
    \fi
  \else
    \ifx\QCBOptA\empty
      {}%
    \else
      \ifx\QCBOptB\empty
        {\QCBOptA}%
      \else
        [\QCBOptB]{\QCBOptA}%
      \fi
    \fi
  \fi
}
\newcount\GRAPHICSTYPE
\GRAPHICSTYPE=\z@
\def\GRAPHICSPS#1{%
 \ifcase\GRAPHICSTYPE
   \special{ps: #1}%
 \or
   \special{language "PS", include "#1"}%
 \fi
}%
\def\graffile#1#2#3#4{%
    \bgroup
	   \@inlabelfalse
       \leavevmode
       \@ifundefined{bbl@deactivate}{\def~{\string~}}{\activesoff}%
        \raise -#4 \BOXTHEFRAME{%
           \hbox to #2{\raise #3\hbox to #2{\null #1\hfil}}}%
    \egroup
}%
\def\draftbox#1#2#3#4{%
 \leavevmode\raise -#4 \hbox{%
  \frame{\rlap{\protect\tiny #1}\hbox to #2%
   {\vrule height#3 width\z@ depth\z@\hfil}%
  }%
 }%
}%
\newcount\@msidraft
\@msidraft=\z@
\let\nographics=\@msidraft
\newif\ifwasdraft
\wasdraftfalse

\def\GRAPHIC#1#2#3#4#5{%
   \ifnum\@msidraft=\@ne\draftbox{#2}{#3}{#4}{#5}%
   \else\graffile{#1}{#3}{#4}{#5}%
   \fi
}
\def\addtoLaTeXparams#1{%
    \edef\LaTeXparams{\LaTeXparams #1}}%
\newif\ifBoxFrame \BoxFramefalse
\newif\ifOverFrame \OverFramefalse
\newif\ifUnderFrame \UnderFramefalse

\def\BOXTHEFRAME#1{%
   \hbox{%
      \ifBoxFrame
         \frame{#1}%
      \else
         {#1}%
      \fi
   }%
}

\def\doFRAMEparams#1{\BoxFramefalse\OverFramefalse\UnderFramefalse\readFRAMEparams#1\end}%
\def\readFRAMEparams#1{%
 \ifx#1\end%
  \let\next=\relax
  \else
  \ifx#1i\dispkind=\z@\fi
  \ifx#1d\dispkind=\@ne\fi
  \ifx#1f\dispkind=\tw@\fi
  \ifx#1t\addtoLaTeXparams{t}\fi
  \ifx#1b\addtoLaTeXparams{b}\fi
  \ifx#1p\addtoLaTeXparams{p}\fi
  \ifx#1h\addtoLaTeXparams{h}\fi
  \ifx#1X\BoxFrametrue\fi
  \ifx#1O\OverFrametrue\fi
  \ifx#1U\UnderFrametrue\fi
  \ifx#1w
    \ifnum\@msidraft=1\wasdrafttrue\else\wasdraftfalse\fi
    \@msidraft=\@ne
  \fi
  \let\next=\readFRAMEparams
  \fi
 \next
 }%
\def\IFRAME#1#2#3#4#5#6{%
      \bgroup
      \let\QCTOptA\empty
      \let\QCTOptB\empty
      \let\QCBOptA\empty
      \let\QCBOptB\empty
      #6%
      \parindent=0pt
      \leftskip=0pt
      \rightskip=0pt
      \setbox0=\hbox{\QCBOptA}%
      \@tempdima=#1\relax
      \ifOverFrame
          \typeout{This is not implemented yet}%
          \show\HELP
      \else
         \ifdim\wd0>\@tempdima
            \advance\@tempdima by \@tempdima
            \ifdim\wd0 >\@tempdima
               \setbox1 =\vbox{%
                  \unskip\hbox to \@tempdima{\hfill\GRAPHIC{#5}{#4}{#1}{#2}{#3}\hfill}%
                  \unskip\hbox to \@tempdima{\parbox[b]{\@tempdima}{\QCBOptA}}%
               }%
               \wd1=\@tempdima
            \else
               \textwidth=\wd0
               \setbox1 =\vbox{%
                 \noindent\hbox to \wd0{\hfill\GRAPHIC{#5}{#4}{#1}{#2}{#3}\hfill}\\%
                 \noindent\hbox{\QCBOptA}%
               }%
               \wd1=\wd0
            \fi
         \else
            \ifdim\wd0>0pt
              \hsize=\@tempdima
              \setbox1=\vbox{%
                \unskip\GRAPHIC{#5}{#4}{#1}{#2}{0pt}%
                \break
                \unskip\hbox to \@tempdima{\hfill \QCBOptA\hfill}%
              }%
              \wd1=\@tempdima
           \else
              \hsize=\@tempdima
              \setbox1=\vbox{%
                \unskip\GRAPHIC{#5}{#4}{#1}{#2}{0pt}%
              }%
              \wd1=\@tempdima
           \fi
         \fi
         \@tempdimb=\ht1
         \advance\@tempdimb by -#2
         \advance\@tempdimb by #3
         \leavevmode
         \raise -\@tempdimb \hbox{\box1}%
      \fi
      \egroup%
}%
\def\DFRAME#1#2#3#4#5{%
  \hfil\break
  \bgroup
     \leftskip\@flushglue
	 \rightskip\@flushglue
	 \parindent\z@
	 \parfillskip\z@skip
     \let\QCTOptA\empty
     \let\QCTOptB\empty
     \let\QCBOptA\empty
     \let\QCBOptB\empty
	 \vbox\bgroup
        \ifOverFrame 
           #5\QCTOptA\par
        \fi
        \GRAPHIC{#4}{#3}{#1}{#2}{\z@}%
        \ifUnderFrame 
           \break#5\QCBOptA
        \fi
	 \egroup
   \egroup
   \break
}%
\def\FFRAME#1#2#3#4#5#6#7{%
  \@ifundefined{floatstyle}
    {
     \begin{figure}[#1]%
    }
    {
	 \ifx#1h
      \begin{figure}[H]%
	 \else
      \begin{figure}[#1]%
	 \fi
	}
  \let\QCTOptA\empty
  \let\QCTOptB\empty
  \let\QCBOptA\empty
  \let\QCBOptB\empty
  \ifOverFrame
    #4
    \ifx\QCTOptA\empty
    \else
      \ifx\QCTOptB\empty
        \caption{\QCTOptA}%
      \else
        \caption[\QCTOptB]{\QCTOptA}%
      \fi
    \fi
    \ifUnderFrame\else
      \label{#5}%
    \fi
  \else
    \UnderFrametrue%
  \fi
  \begin{center}\GRAPHIC{#7}{#6}{#2}{#3}{\z@}\end{center}%
  \ifUnderFrame
    #4
    \ifx\QCBOptA\empty
      \caption{}%
    \else
      \ifx\QCBOptB\empty
        \caption{\QCBOptA}%
      \else
        \caption[\QCBOptB]{\QCBOptA}%
      \fi
    \fi
    \label{#5}%
  \fi
  \end{figure}%
 }%
\newcount\dispkind%

\def\makeactives{
  \catcode`\"=\active
  \catcode`\;=\active
  \catcode`\:=\active
  \catcode`\'=\active
  \catcode`\~=\active
}
\bgroup
   \makeactives
   \gdef\activesoff{%
      \def"{\string"}
      \def;{\string;}
      \def:{\string:}
      \def'{\string'}
      \def~{\string~}
    }
\egroup

\def\FRAME#1#2#3#4#5#6#7#8{%
 \bgroup
 \ifnum\@msidraft=\@ne
   \wasdrafttrue
 \else
   \wasdraftfalse%
 \fi
 \def\LaTeXparams{}%
 \dispkind=\z@
 \def\LaTeXparams{}%
 \doFRAMEparams{#1}%
 \ifnum\dispkind=\z@\IFRAME{#2}{#3}{#4}{#7}{#8}{#5}\else
  \ifnum\dispkind=\@ne\DFRAME{#2}{#3}{#7}{#8}{#5}\else
   \ifnum\dispkind=\tw@
    \edef\@tempa{\noexpand\FFRAME{\LaTeXparams}}%
    \@tempa{#2}{#3}{#5}{#6}{#7}{#8}%
    \fi
   \fi
  \fi
  \ifwasdraft\@msidraft=1\else\@msidraft=0\fi{}%
  \egroup
 }%
\def\TEXUX#1{"texux"}

\long\def\QQQ#1#2{%
     \long\expandafter\def\csname#1\endcsname{#2}}%
\@ifundefined{QTP}{\def\QTP#1{}}{}
\@ifundefined{QEXCLUDE}{\def\QEXCLUDE#1{}}{}
\@ifundefined{Qlb}{}{}
\@ifundefined{Qlt}{}{}
\long\def\QQA#1#2{}%
\def\QTR#1#2{{\csname#1\endcsname #2}}
\def\EXPAND#1[#2]#3{}%
\def\NOEXPAND#1[#2]#3{}%
\def\LaTeXparent#1{}%
\def\ChildStyles#1{}%
\def\ChildDefaults#1{}%
\def\QTagDef#1#2#3{}%

\@ifundefined{correctchoice}{}{}
\@ifundefined{HTML}{\def\HTML#1{\relax}}{}
\@ifundefined{TCIIcon}{\def\TCIIcon#1#2#3#4{\relax}}{}
\if@compatibility
\else
  \providecommand{\UNICODE}[2][]{\protect\rule{.1in}{.1in}}
  \providecommand{\U}[1]{\protect\rule{.1in}{.1in}}
  
\fi

\@ifundefined{lambdabar}{
      
   }{}

\@ifundefined{StyleEditBeginDoc}{}{}
\def\QQfnmark#1{\footnotemark}

\@ifundefined{TCIMAKEINDEX}{}{\makeindex}%
\@ifundefined{abstract}{%
 \def\abstract{%
  \if@twocolumn
   \section*{Abstract (Not appropriate in this style!)}%
   \else \small 
   \begin{center}{\bf Abstract\vspace{-.5em}\vspace{\z@}}\end{center}%
   \quotation 
   \fi
  }%
 }{%
 }%
\@ifundefined{endabstract}{\def\endabstract
  {\if@twocolumn\else\endquotation\fi}}{}%
\@ifundefined{maketitle}{\def\maketitle#1{}}{}%
\@ifundefined{affiliation}{\def\affiliation#1{}}{}%
\@ifundefined{proof}{}{}%
\@ifundefined{endproof}{}{}%
\@ifundefined{newfield}{\def\newfield#1#2{}}{}%
\@ifundefined{chapter}{\def\chapter#1{\par(Chapter head:)#1\par }%
 \newcount\c@chapter}{}%
\@ifundefined{part}{\def\part#1{\par(Part head:)#1\par }}{}%
\@ifundefined{section}{\def\section#1{\par(Section head:)#1\par }}{}%
\@ifundefined{subsection}{\def\subsection#1%
 {\par(Subsection head:)#1\par }}{}%
\@ifundefined{subsubsection}{\def\subsubsection#1%
 {\par(Subsubsection head:)#1\par }}{}%
\@ifundefined{paragraph}{\def\paragraph#1%
 {\par(Subsubsubsection head:)#1\par }}{}%
\@ifundefined{subparagraph}{\def\subparagraph#1%
 {\par(Subsubsubsubsection head:)#1\par }}{}%
\@ifundefined{therefore}{}{}%
\@ifundefined{backepsilon}{}{}%
\@ifundefined{yen}{}{}%
\@ifundefined{registered}{%
   \def\registered{\relax\ifmmode{}\r@gistered
                    \else$\m@th\r@gistered$\fi}%
 \def\r@gistered{^{\ooalign
  {\hfil\raise.07ex\hbox{$\scriptstyle\rm\text{R}$}\hfil\crcr
  \mathhexbox20D}}}}{}%
\@ifundefined{Eth}{}{}%
\@ifundefined{eth}{}{}%
\@ifundefined{Thorn}{}{}%
\@ifundefined{thorn}{}{}%
\@ifundefined{degree}{}{}%
\newdimen\theight
\@ifundefined{Column}{\def\Column{%
 \vadjust{\setbox\z@=\hbox{\scriptsize\quad\quad tcol}%
  \theight=\ht\z@\advance\theight by \dp\z@\advance\theight by \lineskip
  \kern -\theight \vbox to \theight{%
   \rightline{\rlap{\box\z@}}%
   \vss
   }%
  }%
 }}{}%
\@ifundefined{qed}{\def\qed{%
 \ifhmode\unskip\nobreak\fi\ifmmode\ifinner\else\hskip5\p@\fi\fi
 \hbox{\hskip5\p@\vrule width4\p@ height6\p@ depth1.5\p@\hskip\p@}%
 }}{}%
\@ifundefined{cents}{}{}%
\@ifundefined{tciLaplace}{}{}%
\@ifundefined{tciFourier}{}{}%
\@ifundefined{textcurrency}{}{}%
\@ifundefined{texteuro}{}{}%
\@ifundefined{textfranc}{}{}%
\@ifundefined{textlira}{}{}%
\@ifundefined{textpeseta}{}{}%
\@ifundefined{miss}{\def\miss{\hbox{\vrule height2\p@ width 2\p@ depth\z@}}}{}%
\@ifundefined{vvert}{}{}
\@ifundefined{tcol}{\def\tcol#1{{\baselineskip=6\p@ \vcenter{#1}} \Column}}{}%
\@ifundefined{dB}{}{}
\@ifundefined{mB}{}{}
\@ifundefined{nB}{}{}
\@ifundefined{note}{}{}%
\def\newfmtname{LaTeX2e}
\ifx\fmtname\newfmtname
  \DeclareOldFontCommand{\rm}{\normalfont\rmfamily}{\mathrm}
  \DeclareOldFontCommand{\sf}{\normalfont\sffamily}{\mathsf}
  \DeclareOldFontCommand{\tt}{\normalfont\ttfamily}{\mathtt}
  \DeclareOldFontCommand{\bf}{\normalfont\bfseries}{\mathbf}
  \DeclareOldFontCommand{\it}{\normalfont\itshape}{\mathit}
  \DeclareOldFontCommand{\sl}{\normalfont\slshape}{\@nomath\sl}
  \DeclareOldFontCommand{\sc}{\normalfont\scshape}{\@nomath\sc}
\fi

\def\alpha{{\Greekmath 010B}}%
\def\beta{{\Greekmath 010C}}%
\def\gamma{{\Greekmath 010D}}%
\def\delta{{\Greekmath 010E}}%
\def\epsilon{{\Greekmath 010F}}%
\def\zeta{{\Greekmath 0110}}%
\def\eta{{\Greekmath 0111}}%
\def\theta{{\Greekmath 0112}}%
\def\iota{{\Greekmath 0113}}%
\def\kappa{{\Greekmath 0114}}%
\def\lambda{{\Greekmath 0115}}%
\def\mu{{\Greekmath 0116}}%
\def\nu{{\Greekmath 0117}}%
\def\xi{{\Greekmath 0118}}%
\def\pi{{\Greekmath 0119}}%
\def\rho{{\Greekmath 011A}}%
\def\sigma{{\Greekmath 011B}}%
\def\tau{{\Greekmath 011C}}%
\def\upsilon{{\Greekmath 011D}}%
\def\phi{{\Greekmath 011E}}%
\def\chi{{\Greekmath 011F}}%
\def\psi{{\Greekmath 0120}}%
\def\omega{{\Greekmath 0121}}%
\def\varepsilon{{\Greekmath 0122}}%
\def\vartheta{{\Greekmath 0123}}%
\def\varpi{{\Greekmath 0124}}%
\def\varrho{{\Greekmath 0125}}%
\def\varsigma{{\Greekmath 0126}}%
\def\varphi{{\Greekmath 0127}}%

\def\nabla{{\Greekmath 0272}}
\def\FindBoldGroup{%
   {\setbox0=\hbox{$\mathbf{x\global\edef\theboldgroup{\the\mathgroup}}$}}%
}

\def\Greekmath#1#2#3#4{%
    \if@compatibility
        \ifnum\mathgroup=\symbold
           \mathchoice{\mbox{\boldmath$\displaystyle\mathchar"#1#2#3#4$}}%
                      {\mbox{\boldmath$\textstyle\mathchar"#1#2#3#4$}}%
                      {\mbox{\boldmath$\scriptstyle\mathchar"#1#2#3#4$}}%
                      {\mbox{\boldmath$\scriptscriptstyle\mathchar"#1#2#3#4$}}%
        \else
           \mathchar"#1#2#3#4%
        \fi 
    \else 
        \FindBoldGroup
        \ifnum\mathgroup=\theboldgroup 
           \mathchoice{\mbox{\boldmath$\displaystyle\mathchar"#1#2#3#4$}}%
                      {\mbox{\boldmath$\textstyle\mathchar"#1#2#3#4$}}%
                      {\mbox{\boldmath$\scriptstyle\mathchar"#1#2#3#4$}}%
                      {\mbox{\boldmath$\scriptscriptstyle\mathchar"#1#2#3#4$}}%
        \else
           \mathchar"#1#2#3#4%
        \fi     	    
	  \fi}

\newif\ifGreekBold  \GreekBoldfalse
\let\SAVEPBF=\pbf
\def\pbf{\GreekBoldtrue\SAVEPBF}%

\@ifundefined{theorem}{\newtheorem{theorem}{Theorem}}{}
\@ifundefined{lemma}{\newtheorem{lemma}[theorem]{Lemma}}{}
\@ifundefined{corollary}{}{}
\@ifundefined{conjecture}{}{}
\@ifundefined{proposition}{\newtheorem{proposition}[theorem]{Proposition}}{}
\@ifundefined{axiom}{}{}
\@ifundefined{remark}{\newtheorem{remark}{Remark}}{}
\@ifundefined{example}{}{}
\@ifundefined{exercise}{}{}
\@ifundefined{definition}{\newtheorem{definition}{Definition}}{}

\@ifundefined{mathletters}{%
  \newcounter{equationnumber}  
  \def\mathletters{%
     \addtocounter{equation}{1}
     \edef\@currentlabel{\theequation}%
     \setcounter{equationnumber}{\c@equation}
     \setcounter{equation}{0}%
     \edef\theequation{\@currentlabel\noexpand\alph{equation}}%
  }
  
}{}

\@ifundefined{BibTeX}{%
    \def\BibTeX{{\rm B\kern-.05em{\sc i\kern-.025em b}\kern-.08em
                 T\kern-.1667em\lower.7ex\hbox{E}\kern-.125emX}}}{}%
\@ifundefined{AmS}%
    {\def\AmS{{\protect\usefont{OMS}{cmsy}{m}{n}%
                A\kern-.1667em\lower.5ex\hbox{M}\kern-.125emS}}}{}%
\@ifundefined{AmSTeX}{}{}%
\def\@@eqncr{\let\@tempa\relax
    \ifcase\@eqcnt \def\@tempa{& & &}\or \def\@tempa{& &}%
      \else \def\@tempa{&}\fi
     \@tempa
     \if@eqnsw
        \iftag@
           \@taggnum
        \else
           \@eqnnum\stepcounter{equation}%
        \fi
     \fi
     \global\tag@false
     \global\@eqnswtrue
     \global\@eqcnt\z@\cr}

\def\TCItag{\@ifnextchar*{\@TCItagstar}{\@TCItag}}
\def\@TCItag#1{%
    \global\tag@true
    \global\def\@taggnum{(#1)}}
\def\@TCItagstar*#1{%
    \global\tag@true
    \global\def\@taggnum{#1}}
\RequirePackage{amsmath}
\makeatother

\begin{document}

\title{Equilibria for a gyrostat in newtonian interaction with two rigid bodies}
\author{J. A. Vera\\Universidad Polit\'{e}cnica de Cartagena\\Departamento de Matem\'{a}tica Aplicada y Estad\'{\i}stica \\Paseo Alfonso XIII, 52, 30203 Cartagena (Murcia), Spain\\e-mail: juanantonio.vera@upct.es}
\date{}
\maketitle

\begin{abstract}
In this paper the non-canonical Hamiltonian dynamics of a gyrostat
in the three body problem will be examined. By means of
geometric-mechanics methods some relative equilibria of the
dynamics of a gyrostat in Newtonian interaction with two rigid
bodies will be studied. Taking advantage of the results obtained
in previous papers, working on the reduced problem, the
bifurcations of these relative equilibria will be studied. The
instability of Eulerian relative equilibria if the gyrostat is
close to a sphere is proven. Necessary and sufficient conditions
will be provided for lineal stability of Lagrangian relative
equilibria if the gyrostat is close to a sphere. The analysis is
done in vectorial form avoiding the use of canonical variables and
the tedious expressions associated with them. In this way, the
classic results on equilibria of the three-body problem, many of
them obtained by other authors who had used of more classic
techniques, are generalized.


\end{abstract}

\DeclareGraphicsExtensions{.jpg}

\pagestyle{fancy} \lhead{\small J. A. Vera} \rhead{\small Equilibria for a
gyrostat in newtonian interaction with two rigid bodies}


\section{Introduction}

In the last few years a lot of papers about the problem of roto-translational
motion of celestial bodies have appeared. These show a new interest in the
study of configurations of relative equilibria by differential geometry
methods instead of more classical ones.

In the present paper, references will be made to Wang et al. \cite{9},
concernig the problem of a rigid body in a central Newtonian field, and
Maciejewski \cite{3}, concerning the problem of two rigid bodies in mutual
Newtonian attraction. Let us remember that a gyrostat is a mechanical system
$\mathit{S}$ composed of a rigid body $\mathit{S^{\prime}}$ and other bodies
$\mathit{S^{\prime\prime}}$ (deformable or rigid), connected to it in such a
way that their relative motion with respect to its rigid part do not change
the distribution of mass of the total system $\mathit{S}$ (Leimanis,
\cite{2}). These papers have been generalized by Mond\'{e}jar and Vigueras
\cite{4} to the case of two gyrostats in mutual Newtonian attraction.

With regards to the problem of three rigid bodies, references will have to be
made to Vidiakin \cite{8} and Duboshin \cite{1} who proved the existence of
Euler and Lagrange configurations of equilibria when the bodies possess
symmetries. Later Zhuravlev and Petrutskii \cite{10} reviewed the results
until 1990.

Vera \cite{5} and Vera and Vigueras (\cite{6},\cite{7}) have studied the
non-canonical Hamiltonian dynamics of $n+1$ bodies in Newtonian attraction, in
which $n$ of them are rigid bodies with a spherical distribution of masses or
material points and the other one is a triaxial gyrostat. Working in the
reduced problem, global considerations about the conditions for relative
equilibria were made. Finally, in an approximated model of the dynamics, a
study of some relative equilibria of a gyrostat in Newtonian interaction with
three rigid bodies is carried out.

In this paper, the existence and stability of some periodic solutions of the
dynamics of a gyrostat in Newtonian interaction with two rigid bodies will be
investigated. In a first approach to the problem, the attitude dynamics of
both rigid bodies are assumed to be the same as the dynamics of a rigid body
in torque free motion. The two rigid bodies have revolution symmetry about the
third axis of inertia. On the other hand, the Newtonian interaction of the
gyrostat between both rigid bodies is simplified supposing the gyrostat a
material particle. The Newtonian interaction between the two rigid bodies are
the same of two material particles.

With these hypotheses, according to Vera and Vigueras \cite{7}, working in the
double reduced space of configuration of the problem, the equations of motion
and those which determine the relative equilibria will be derived.

Two families of relative equilibria, \emph{Eulerian} and \emph{Lagrangian}
will be investigated. The Eulerian relative equilibria are completely
determinated by a polynomial equation of degree nine. The bifurcations of
these equilibria are made. The bifurcations of the Lagrangian relative
equilibria are completely investigated and expressions for the Lagrangian
relative equilibria when both solids are close to spheres are provided.

We should notice that the studied system, has potential interest both in
astrodynamics (dealing with spacecrafts) as well as in the understanding of
the evolution of planetary systems recently found (and more to appear), where
some of the planets may be modeled like a rigid body rather than a rigid body.
In fact, the equilibria reported might well be compared with the ones taken
for the `parking areas' of the space missions (GENESIS, SOHO, DARWIN, etc)
around the Eulerian points of the Sun-Earth and the Earth-Moon systems. Some
interesting numeric results have been calculated in the Appendix B.

The methods used in this work can be applied to similar problems. The study of
the nonlinear stability of the relative equilibria obtained here is the
logical continuation of this work.


\section{Poisson dynamics}

Let $S_{0}$ be a gyrostat of mass $m_{0},$ and $S_{1},$ $S_{2}$ two
symmetrical rigid bodies of masses $m_{1}$ and $m_{2}$ respectively;
$\mathfrak{I}=\{\mathbf{O},\mathbf{u}_{1},\mathbf{u}_{2},\mathbf{u}_{3}\}$ an
inertial reference frame$;\ \mathfrak{J}=\{\mathbf{C}_{0},\mathbf{b}%
_{1},\mathbf{b}_{2},\mathbf{b}_{3}\}$ a body frame fixed at the center of mass
$\mathbf{C}_{0}$ of $S_{0},$ (see Vera and Vigueras \cite{7} for details).

The following notation is used%
\[%
\begin{array}
[c]{c}%
M_{2}=m_{1}+m_{2},\quad M_{1}=m_{1}+m_{2}+m_{0},\quad g_{1}=\dfrac{m_{1}m_{2}%
}{M_{2}},\quad g_{2}=\dfrac{m_{0}M_{2}}{M_{1}}%
\end{array}
\]

Where for $\mathbf{u}$, $\mathbf{v}\in\mathbb{R}^{3}$, $\mathbf{u\cdot v}$ is
the dot product; $\mid\mathbf{u}\mid$ is the Euclidean norm of the vector
$\mathbf{u}$; $\mathbf{u\times v}$ is the cross product; $\mathbf{I}%
_{\mathbb{R}^{3}}$ is the identity matrix; and $\mathbf{0}$ is the zero matrix
of order three. Consider $\mathbb{I}_{0}=diag(A_{0},$ $A_{0},$ $C_{0})$ the
diagonal tensor of inertia of the gyrostat, and $\mathbb{I}_{i}=diag(A_{i},$
$A_{i},$ $C_{i})$ the diagonal tensors of inertia for the rigid bodies
$S_{i},$ $i=1,2.$ The generic expression $\mathbf{z=}(\mathbf{\mathbf{\Pi}%
}_{1}\mathbf{,}$ $\mathbf{\mathbf{\Pi}}_{2}\mathbf{,}$ $\mathbf{\Pi}%
_{0}\mathbf{,}$ $\mathbf{\lambda},$ $\mathbf{p}_{\mathbf{\lambda}},$
$\mathbf{\mu},$ $\mathbf{p}_{\mu})\in\mathbb{R}^{21}$ is an vector of the
twice reduced problem obtained by appling the symmetries of the system.
$\mathbf{\Pi}_{0}=\mathbb{I}_{0}\mathbf{\Omega}_{0}+\mathbf{l}_{r}$ is the
total rotational angular momentum vector of the gyrostat in the body frame,
which is attached to its rigid part and whose axes have the direction of the
principal axes of inertia of $S_{0};$ the vector $\mathbf{l}_{r}=(0,0,l)$ is
the constant gyrostatic momentum and $\mathbf{\Pi}_{i}=\mathbb{I}%
_{i}\mathbf{\Omega}_{i}$ $(i=1,2)$ are the total rotational angular momentum
vectors for the two rigid bodies. The elements $\mathbf{\lambda}$,
$\mathbf{\mu}$, $\mathbf{p}_{\mathbf{\lambda}}$ and $\mathbf{p}_{\mu}$ are
respectively the barycentric coordinates and the linear momenta expressed in
the body frame $\mathfrak{J}$.

Following the results of Vera and Vigueras \cite{7}, according to the
hypotheses formulated in the introduction of this paper, a good approximation
to the potential of the system is expressed by the following expression%
\[
\mathcal{V}=\mathcal{V}_{1}+\mathcal{V}_{2}%
\]
where%
\begin{align*}
\mathcal{V}_{1}  &  =-\left(  \dfrac{Gm_{1}m_{2}}{\mid\mathbf{\lambda}\mid
}+\dfrac{Gm_{1}m_{0}}{\mid\mathbf{\mu-}\frac{m_{2}}{M_{2}}\mathbf{\lambda}%
\mid}+\dfrac{Gm_{2}m_{0}}{\mid\mathbf{\mu+}\frac{m_{1}}{M_{2}}\mathbf{\lambda
}\mid}\right)  \medskip\\
\mathcal{V}_{2}  &  =-\dfrac{1}{2}\left(  \dfrac{Gm_{0}\alpha_{1}}%
{\mid\mathbf{\mu-}\frac{m_{2}}{M_{2}}\mathbf{\lambda}\mid^{3}}+\dfrac
{Gm_{0}\alpha_{2}}{\mid\mathbf{\mu+}\frac{m_{1}}{M_{2}}\mathbf{\lambda}%
\mid^{3}}-\dfrac{3Gm_{0}f_{1}}{\mid\mathbf{\mu-}\frac{m_{2}}{M_{2}%
}\mathbf{\lambda}\mid^{5}}-\dfrac{3Gm_{0}f_{2}}{\mid\mathbf{\mu+}\frac{m_{1}%
}{M_{2}}\mathbf{\lambda}\mid^{5}}\right)
\end{align*}
and%
\begin{align*}
\alpha_{1}  &  =2A_{1}+C_{1},\quad\alpha_{2}=2A_{2}+C_{2}\medskip\\
f_{1}(\mathbf{\lambda},\mathbf{\mu})  &  =\mathbf{\mu}\cdot\mathbb{I}%
_{1}\mathbf{\mu}-\dfrac{2m_{2}}{M_{2}}\mathbf{\lambda}\cdot\mathbb{I}%
_{1}\mathbf{\mu}+\left(  \dfrac{m_{2}}{M_{2}}\right)  ^{2}\mathbf{\lambda
}\cdot\mathbb{I}\mathbf{\lambda}\medskip\\
f_{2}(\mathbf{\lambda},\mathbf{\mu})  &  =\mathbf{\mu}\cdot\mathbb{I}%
_{2}\mathbf{\mu}+\dfrac{2m_{1}}{M_{2}}\mathbf{\lambda}\cdot\mathbb{I}%
_{2}\mathbf{\mu}+\left(  \dfrac{m_{1}}{M_{2}}\right)  ^{2}\mathbf{\lambda
}\cdot\mathbb{I}\mathbf{\lambda}%
\end{align*}

The Hamiltonian function of the system adopts the form%
\[
\mathcal{H}(\mathbf{z})=\dfrac{\mid\mathbf{p}_{\mathbf{\lambda}}\mid^{2}%
}{2g_{1}}+\dfrac{\mid\mathbf{p}_{\mathbf{\mu}}\mid^{2}}{2g_{2}}+\dfrac{1}%
{2}\Pi_{0}\mathbb{I}_{0}^{-1}\Pi_{0}-\mathbf{l}_{r}\cdot\mathbb{I}_{0}^{-1}%
\Pi\medskip+\dfrac{1}{2}\Pi_{1}\mathbb{I}_{1}^{-1}\Pi_{1}+\dfrac{1}{2}\Pi
_{2}\mathbb{I}_{2}^{-1}\Pi_{2}+\mathcal{V}(\mathbf{\lambda,\mu})
\]

Let $(\mathbf{M},\{\ ,\ \}$,$\mathcal{H}$$)$ with $\mathbf{M}=\mathbb{R}^{21}
$ be the Poisson manifold, where $\{\ ,\ \}$

is the Poisson brackets defined by means of the Poisson tensor%
\[
\mathbf{B(z)}=\left(
\begin{array}
[c]{ccccccc}%
\widehat{\mathbf{\Pi}_{1}} & \mathbf{0} & \mathbf{0} & \mathbf{0} & \mathbf{0}
& \mathbf{0} & \mathbf{0}\\
\mathbf{0} & \widehat{\mathbf{\Pi}_{2}} & \mathbf{0} & \mathbf{0} & \mathbf{0}
& \mathbf{0} & \mathbf{0}\\
\mathbf{0} & \mathbf{0} & \widehat{\mathbf{\Pi}_{0}} & \widehat
{\mathbf{\lambda}} & \widehat{\mathbf{p}_{\mathbf{\lambda}}} & \widehat
{\mathbf{\mu}} & \widehat{\mathbf{p}_{\mu}}\\
\mathbf{0} & \mathbf{0} & \widehat{\mathbf{\lambda}_{1}} & \mathbf{0} &
\mathbf{I}_{\mathbb{R}^{3}} & \mathbf{0} & \mathbf{0}\\
\mathbf{0} & \mathbf{0} & \widehat{\mathbf{p}_{\mathbf{\lambda}}} &
-\mathbf{I}_{\mathbb{R}^{3}} & \mathbf{0} & \mathbf{0} & \mathbf{0}\\
\mathbf{0} & \mathbf{0} & \widehat{\mathbf{\mu}} & \mathbf{0} & \mathbf{0} &
\mathbf{0} & \mathbf{I}_{\mathbb{R}^{3}}\\
\mathbf{0} & \mathbf{0} & \widehat{\mathbf{p}_{\mu}} & \mathbf{0} & \mathbf{0}
& -\mathbf{I}_{\mathbb{R}^{3}} & \mathbf{0}%
\end{array}
\right)
\]

In $\mathbf{B(z)}$, $\widehat{\mathbf{v}}$ is considered to be the image of
the vector $v\in\mathbb{R}^{3}$ by the standard isomorphism between the Lie
Algebras $\mathbb{R}^{3}$ and $\mathfrak{so(3)}$, i.e.%
\[%
\begin{array}
[c]{c}%
\widehat{\mathbf{v}}=\left(
\begin{array}
[c]{ccc}%
0 & -v_{3} & v_{2}\\
v_{3} & 0 & -v_{1}\\
-v_{2} & v_{1} & 0
\end{array}
\right)
\end{array}
\]

The equations of the motion are given by the below expression%
\[%
\begin{array}
[c]{c}%
\dfrac{d\mathbf{z}}{dt}=\{\mathbf{z},\mathcal{H}(\mathbf{z})\}=\mathbf{B}%
(\mathbf{z)\nabla}_{\mathbf{z}}\mathcal{H}(\mathbf{z})
\end{array}
\]
with $\mathbf{\nabla}_{\mathbf{z}}\mathcal{V}$ being the gradient of
$\mathcal{V}$ with respect to an arbitrary vector $\mathbf{z}$.

Calculating $\{\mathbf{z},\mathcal{H}(\mathbf{z})\}$, the below group of
vectorial equations of the motion can be written as%
\begin{equation}%
\begin{array}
[c]{l}%
\dfrac{d\mathbf{\Pi}_{0}}{dt}=\mathbf{\Pi}_{0}\times\mathbf{\mathbf{\Omega}%
}_{0}\,\mathbf{+\,\lambda\times\nabla}_{\mathbf{\lambda}}\mathcal{V}%
+\mathbf{\mu\times\nabla}_{\mathbf{\mu}}\mathcal{V}\medskip\\
\dfrac{d\mathbf{\lambda}}{dt}=\dfrac{\mathbf{p}_{\mathbf{\lambda}}}{g_{1}%
}+\mathbf{\lambda\times\Omega}_{0}\mathbf{,}\quad\dfrac{d\mathbf{p}%
_{\mathbf{\lambda}}}{dt}=\mathbf{p}_{\mathbf{\lambda}}\times\mathbf{\Omega
}_{0}-\mathbf{\nabla}_{\mathbf{\lambda}}\mathcal{V}\medskip\\
\dfrac{d\mathbf{\mu}}{dt}=\dfrac{\mathbf{p}_{\mathbf{\mu}}}{g_{2}}%
+\mathbf{\mu\times\Omega}_{0}\mathbf{,}\quad\dfrac{d\mathbf{p}_{\mathbf{\mu}}%
}{dt}=\mathbf{p}_{\mathbf{\mu}}\times\mathbf{\Omega}_{0}-\mathbf{\nabla
}_{\mathbf{\mu}}\mathcal{V}\medskip\\
\dfrac{d\mathbf{\Pi}_{1}}{dt}=\mathbf{\Pi}_{1}\times\mathbf{\mathbf{\Omega}%
}_{1}\mathbf{,}\quad\dfrac{d\mathbf{\Pi}_{2}}{dt}=\mathbf{\Pi}_{2}%
\times\mathbf{\mathbf{\Omega}}_{2}%
\end{array}
\label{EcuacHamilt}%
\end{equation}

Important elements of $\mathbf{B(z)}$ are the associated Casimir functions.
The vector%
\[%
\begin{array}
[c]{c}%
\mathbf{L}_{0}=\mathbf{\Pi}_{0}\,\mathbf{+\,\lambda}\times\mathbf{p}%
_{\mathbf{\lambda}}+\mathbf{\mu}\times\mathbf{p}_{\mu}%
\end{array}
\]
is a part of the total angular momentum $\mathbf{L}$ given by%
\[%
\begin{array}
[c]{c}%
\mathbf{L}=\mathbf{\Pi}_{2}+\mathbf{\Pi}_{1}+\mathbf{L}_{0}%
\end{array}
\]
Then the below result can be concluded.

\begin{proposition}
If $\varphi_{i},$ $(i=0,$ $1,$ $2)$ are real smooth functions, then
$\varphi_{0}(\tfrac{\mid\mathbf{L}_{0}\mathbf{\mid}^{2}}{2}),$ $\varphi
_{i}(\tfrac{\mid\mathbf{\Pi}_{i}\mathbf{\mid}^{2}}{2}),$ $(i=1,$ $2)$ are
Casimir functions of the Poisson tensor $\mathbf{B}(\mathbf{z)}$. Furthermore,
$Ker\mathbf{B}(\mathbf{z})=<\mathbf{\nabla}_{\mathbf{z}}\varphi_{0},$
$\mathbf{\nabla}_{\mathbf{z}}\varphi_{1},$ $\mathbf{\nabla}_{\mathbf{z}%
}\varphi_{2}>$. We also have $\tfrac{d\mathbf{L}}{dt}=\mathbf{0}$, which means
the total angular momentum vector remains constant. If $\mathbf{\Pi}%
_{0}=(\mathbf{\pi}_{0}^{1},$ $\mathbf{\pi}_{0}^{2},$ $\mathbf{\pi}_{0}^{3}),$
then $\mathbf{\pi}_{0}^{3}$ is an integral of the motion.
\end{proposition}

\section{Relative Equilibria}

If $\mathbf{z}_{e}=(\mathbf{\Pi}_{2}^{e},$ $\mathbf{\Pi}_{1}^{e},$
$\mathbf{\Pi}_{0}^{e},$ $\mathbf{\lambda}^{e},$ $\mathbf{p}_{\mathbf{\lambda}%
}^{e},$ $\mathbf{\mu}^{e},$ $\mathbf{p}_{\mu}^{e})$ is a generic relative
equilibrium, the below vectorial equations are verified%
\begin{gather}
\mathbf{\Pi}_{0}^{e}\times\mathbf{\mathbf{\Omega}}_{0}^{e}\mathbf{\,+\,\lambda
}^{e}\times(\mathbf{\nabla}_{\mathbf{\lambda}}\mathcal{V})_{e}+\mathbf{\mu
}\times(\mathbf{\nabla}_{\mathbf{\mu}}\mathcal{V})_{e}=\mathbf{0}%
\medskip\label{equilibriosn=2}\\
\dfrac{\mathbf{p}_{\mathbf{\lambda}}^{e}}{g_{1}}+\mathbf{\lambda}^{e}%
\times\mathbf{\Omega}_{0}^{e}=\mathbf{0,}\quad\mathbf{p}_{\mathbf{\lambda}%
}^{e}\times\mathbf{\Omega}_{0}^{e}=(\mathbf{\nabla}_{\mathbf{\lambda}%
}\mathcal{V})_{e}\medskip\nonumber\\
\dfrac{\mathbf{p}_{\mathbf{\mu}}^{e}}{g_{2}}+\mathbf{\mu}^{e}\times
\mathbf{\Omega}_{0}^{e}=\mathbf{0,}\quad\mathbf{p}_{\mathbf{\mu}}^{e}%
\times\mathbf{\Omega}_{0}^{e}=(\mathbf{\nabla}_{\mathbf{\mu}}\mathcal{V}%
)_{e}\medskip\nonumber\\
\mathbf{\Pi}_{1}^{e}\times\mathbf{\mathbf{\Omega}}_{1}^{e}=\mathbf{0,}%
\quad\mathbf{\Pi}_{2}^{e}\times\mathbf{\mathbf{\Omega}}_{2}^{e}=\mathbf{0}%
\end{gather}
where $(\mathbf{\nabla}_{\mathbf{\lambda}}\mathcal{V})_{e}$ and
$(\mathbf{\nabla}_{\mathbf{\mu}}\mathcal{V})_{e}$ are the values of
$\mathbf{\nabla}_{\mathbf{\lambda}}\mathcal{V}$ and $\mathbf{\nabla
}_{\mathbf{\mu}}\mathcal{V} $ in $\mathbf{z}_{e}.$

According to the relationships provided by Vera and Vigueras \cite{7}, the
following results are obtained.

\begin{lemma}
Whenever $\mathbf{z}_{e}=(\mathbf{\Pi}_{2}^{e},$ $\mathbf{\Pi}_{1}^{e},$
$\mathbf{\Pi}_{0}^{e},$ $\mathbf{\lambda}^{e},$ $\mathbf{p}_{\mathbf{\lambda}%
}^{e},$ $\mathbf{\mu}^{e},$ $\mathbf{p}_{\mu}^{e})$ is a relative equilibrium,
the below relationships are verified%
\begin{align*}
&  \mid\mathbf{\mathbf{\Omega}}_{0}^{e}\mid^{2}\mid\mathbf{\lambda}^{e}%
\mid^{2}-(\mathbf{\lambda}^{e}\cdot\mathbf{\mathbf{\Omega}}_{0}^{e}%
)^{2}=\dfrac{1}{g_{1}}(\mathbf{\lambda}^{e}\cdot(\mathbf{\nabla}%
_{\mathbf{\lambda}}\mathcal{V})_{e})\medskip\\
&  \mid\mathbf{\mathbf{\Omega}}_{0}^{e}\mid^{2}\mid\mathbf{\mu}^{e}\mid
^{2}-(\mathbf{\mu}^{e}\cdot\mathbf{\mathbf{\Omega}}_{0}^{e})^{2}=\dfrac
{1}{g_{2}}(\mathbf{\mu}^{e}\cdot(\mathbf{\nabla}_{\mathbf{\mu}}\mathcal{V}%
)_{e})
\end{align*}

\end{lemma}

The previous two identities will be used to obtain necessary conditions for
the existence of relative equilibria.

Certain relative equilibria will be studied assuming that vectors
$\mathbf{\mathbf{\Omega}}_{0}^{e},$ $\mathbf{\lambda}^{e}$ and $\mathbf{\mu
}^{e}$ satisfy special geometric properties.

\begin{definition}
$\mathbf{z}_{e}$ is said to be an \emph{Eulerian relative equilibrium} when,
$\mathbf{\lambda}^{e}$ and $\mathbf{\mu}^{e}$ are proportional and
$\mathbf{\Omega}_{e}$ is perpendicular to the straight line that is generate.
\end{definition}

\begin{definition}
$\mathbf{z}_{e}$ is said to be a \emph{Lagrangian relative equilibrium} if
$\mathbf{\lambda}^{e}$ and $\mathbf{\mu}^{e}$ are not proportional and
$\mathbf{\Omega}_{e}$ is perpendicular to the plane that both of them generate.
\end{definition}

From the above definitions, the following property is deduced.

\begin{proposition}
In an Eulerian or Lagrangian relative equilibrium, momenta are not exercised
over the gyrostat.
\end{proposition}

Next, necessary and sufficient conditions for the existence of Eulerian and
Lagrangian relative equilibria will we obtained.

\section{Eulerian relative equilibria}

According to the relative position of the gyrostat $S_{0}$ with respect to
$S_{1}$ and $S_{2},$ there are three possible equilibrium configurations: a)
$S_{2}S_{1}S_{0}$, b) $S_{2}S_{0}S_{1}$ and c) $S_{0}S_{2}S_{1}$.

\subsection{Necessary conditions of existence}

The following lemma is a direct consequence of the geometry of the problem.

\begin{lemma}
If $\mathbf{z}_{e}=(\mathbf{\Pi}_{2}^{e},$ $\mathbf{\Pi}_{1}^{e},$
$\mathbf{\Pi}_{0}^{e},$ $\mathbf{\lambda}^{e},$ $\mathbf{p}_{\mathbf{\lambda}%
}^{e},$ $\mathbf{\mu}^{e},$ $\mathbf{p}_{\mu}^{e})$ is a relative equilibrium
of Euler type, then for the configuration $S_{0}S_{2}S_{1}$%
\[
\mid\mathbf{\mu}^{e}\mathbf{-}\dfrac{m_{2}}{M_{2}}\mathbf{\lambda}^{e}%
\mid\,=\,\mid\mathbf{\lambda}^{e}\mid+\mid\mathbf{\mu}^{e}\mathbf{+}%
\dfrac{m_{1}}{M_{2}}\mathbf{\lambda}^{e}\mid
\]
In a similar way, for the configuration $S_{2}S_{0}S_{1}$%
\[
\mid\mathbf{\lambda}^{e}\mid\,=\,\mid\mathbf{\mu}^{e}\mathbf{-}\dfrac{m_{2}%
}{M_{2}}\mathbf{\lambda}^{e}\mid+\mid\mathbf{\mu}^{e}\mathbf{+}\dfrac{m_{1}%
}{M_{2}}\mathbf{\lambda}^{e}\mid
\]
Finally, for the configuration $S_{2}S_{1}S_{0}$%
\[
\mid\mathbf{\mu}^{e}\mathbf{+}\dfrac{m_{1}}{M_{2}}\mathbf{\lambda}^{e}%
\mid\,=\,\mid\mathbf{\mu}^{e}\mathbf{-}\dfrac{m_{2}}{M_{2}}\mathbf{\lambda
}^{e}\mid+\mid\mathbf{\lambda}^{e}\mid
\]

\end{lemma}

If $\mathbf{z}_{e}$ is an Eulerian relative equilibrium, them%
\begin{align*}
g_{1}  &  \mid\mathbf{\Omega}_{e}^{0}\mid^{2}\mid\mathbf{\lambda}^{e}\mid
^{2}\,=\mathbf{\lambda}^{e}\cdot(\mathbf{\nabla}_{\mathbf{\lambda}}%
\mathcal{V})_{e}\\
& \\
g_{2}  &  \mid\mathbf{\Omega}_{e}^{0}\mid^{2}\mid\mathbf{\mu}^{e}\mid
^{2}\,=\mathbf{\mu}^{e}\cdot(\mathbf{\nabla}_{\mathbf{\mu}}\mathcal{V})_{e}%
\end{align*}

with
\[
\mathbf{\mu}^{e}\mathbf{-}\dfrac{m_{2}}{M_{2}}\mathbf{\lambda}^{e}%
=\rho\mathbf{\lambda}^{e},\quad\mathbf{\mu}^{e}\mathbf{+}\dfrac{m_{1}}{M_{2}%
}\mathbf{\lambda}^{e}\mathbf{=(}1+\rho)\mathbf{\lambda}^{e},\quad\mathbf{\mu
}^{e}=\dfrac{\left(  (1+\rho)m_{2}+\rho m_{1}\right)  }{M_{2}}\mathbf{\lambda
}^{e}\medskip
\]
where $\rho\in(-\infty,-1)$ in the configuration $a)$, $\rho\in(-1,0)$ in the
configuration $b)$ and $\rho\in(0,+\infty)$ in the configuration $c)$.
Moreover the following expressions are possible to be obtained%
\[
(\mathbf{\nabla}_{\mathbf{\lambda}}\mathcal{V})_{e}=h_{1}(\rho)\mathbf{\lambda
}^{e},\quad(\mathbf{\nabla}_{\mathbf{\mu}}\mathcal{V})_{e}=h_{2}%
(\rho)\mathbf{\lambda}^{e}\medskip
\]
with%
\begin{align}
h_{1}(\rho)  &  =\dfrac{Gm_{1}m_{2}}{\mid\mathbf{\lambda}^{e}\mid^{3}}%
+\dfrac{Gm_{0}m_{1}sgn(1+\rho)}{M_{2}\mid\mathbf{\lambda}^{e}\mid^{3}}\left(
\dfrac{m_{2}}{(1+\rho)^{2}}+\dfrac{\beta_{2}}{(1+\rho)^{4}\mid\mathbf{\lambda
}^{e}\mid^{2}}\right)  -\medskip\label{h1}\\
&  \dfrac{Gm_{0}m_{2}sgn(\rho)}{M_{2}\mid\mathbf{\lambda}^{e}\mid^{3}}\left(
\dfrac{m_{1}}{\rho^{2}}+\dfrac{\beta_{1}}{\rho^{4}\mid\mathbf{\lambda}^{e}%
\mid^{2}}\right)
\end{align}
and%
\begin{align}
h_{2}(\rho)  &  =\dfrac{Gm_{0}sgn(1+\rho)}{\mid\mathbf{\lambda}^{e}\mid^{3}%
}\left(  \dfrac{m_{2}}{(1+\rho)^{2}}+\dfrac{\beta_{2}}{\mid\mathbf{\lambda
}^{e}\mid^{2}(1+\rho)^{4}}\right)  +\medskip\label{h2}\\
&  \dfrac{Gm_{0}sgn(\rho)}{\mid\mathbf{\lambda}^{e}\mid^{3}}\left(
\dfrac{m_{1}}{\rho^{2}}+\dfrac{\beta_{1}}{\mid\mathbf{\lambda}^{e}\mid^{2}%
\rho^{4}}\right)
\end{align}
where $\beta_{1}=3(C_{1}-A_{1})/2,$ $\beta_{2}=3(C_{2}-A_{2})/2$ and
$sgn(x)=\left\{
\begin{array}
[c]{c}%
1\text{ if }x>0\\
-1\text{ if }x\leq0\text{ }%
\end{array}
\right.  $.

Now, from the identities%
\begin{align*}
\mathbf{\lambda}^{e}\cdot(\mathbf{\nabla}_{\mathbf{\lambda}}\mathcal{V})_{e}
&  =\mid\mathbf{\lambda}^{e}\mid^{2}h_{1}(\rho)\medskip\\
\mathbf{\mu}^{e}\cdot(\mathbf{\nabla}_{\mathbf{\mu}}\mathcal{V})_{e}  &
=\dfrac{\left(  (1+\rho)m_{2}+\rho m_{1}\right)  }{M_{2}}\mid\mathbf{\lambda
}^{e}\mid^{2}h_{2}(\rho)
\end{align*}
these equations are deduced%
\begin{align*}
&  \mid\mathbf{\mathbf{\Omega}}_{0}^{e}\mid^{2}=\dfrac{m_{1}+m_{2}}{m_{1}%
m_{2}}h_{1}(\rho)\medskip\\
&  \mid\mathbf{\mathbf{\Omega}}_{0}^{e}\mid^{2}=\dfrac{m_{0}+m_{1}+m_{2}%
}{m_{0}\left(  (1+\rho)m_{2}+\rho m_{1}\right)  }h_{2}(\rho)
\end{align*}
Then for an Eulerian relative equilibrium, $\rho$ must be a real root of the
below equation%
\begin{equation}
m_{0}(m_{1}+m_{2})\left(  (1+\rho)m_{1}+\rho m_{2}\right)  h_{1}(\rho
)=m_{1}m_{2}(m_{0}+m_{1}+m_{2})h_{2}(\rho) \label{EcEuler1}%
\end{equation}

The following proposition summarizes all these results.

\begin{proposition}
If $\mathbf{z}_{e}=(\mathbf{\Pi}_{2}^{e},$ $\mathbf{\Pi}_{1}^{e},$
$\mathbf{\Pi}_{0}^{e},$ $\mathbf{\lambda}^{e},$ $\mathbf{p}_{\mathbf{\lambda}%
}^{e},$ $\mathbf{\mu}^{e},$ $\mathbf{p}_{\mu}^{e})$ is an Eulerian relative
equilibrium, the equation (\ref{EcEuler1}) has, at least, a real root in which
the functions $h_{1}(\rho)$ and $h_{2}(\rho)$ are given by the expressions
(\ref{h1},\ref{h2}). The modulus of the angular velocity of the gyrostat is
\[
\mid\mathbf{\mathbf{\Omega}}_{0}^{e}\mid^{2}=\dfrac{m_{1}+m_{2}}{m_{1}m_{2}%
}h_{1}(\rho)
\]

\begin{remark}
When $\mid\mathbf{\lambda}_{e}\mid$ has a fixed value, if an Eulerian relative
equilibrium exists, the equation (\ref{EcEuler1}) has real solutions. The
number of real roots of the equation (\ref{EcEuler1}) will depend, obviously,
on the parameters which exist in the system.
\end{remark}
\end{proposition}

\subsection{Sufficient conditions of existence}

The below proposition indicates how to calculate solutions of eq.
(\ref{equilibriosn=2}).

\begin{proposition}
When $\mid\mathbf{\lambda}^{e}\mid$ has a fixed value, let $\rho$ be a
solution of the equation (\ref{EcEuler1}), where the functions $h_{1}(\rho)$
and $h_{2}(\rho)$ are given by the expressions (\ref{h1},\ref{h2}). Then
$\mathbf{z}_{e}=(\mathbf{\Pi}_{2}^{e},$ $\mathbf{\Pi}_{1}^{e},$ $\mathbf{\Pi
}_{0}^{e},$ $\mathbf{\lambda}^{e},$ $\mathbf{p}_{\mathbf{\lambda}}^{e}, $
$\mathbf{\mu}^{e},$ $\mathbf{p}_{\mu}^{e})$ given by%
\begin{gather*}
\mathbf{\lambda}^{e}=(\lambda^{e},0,0),\quad\mathbf{\mu}^{e}=(\mu
^{e},0,0),\quad\mathbf{\mathbf{\Omega}}_{0}^{e}=(0,0,\omega_{0}^{e})\medskip\\
\mathbf{p}_{\mathbf{\lambda}}^{e}=(0,g_{1}\omega_{0}^{e}\lambda^{e}%
,0),\quad\mathbf{p}_{\mathbf{\mu}}^{e}=(0,g_{2}\omega_{0}^{e}\mu^{e}%
,0),\quad\mathbf{\Pi}_{0}^{e}=(0,0,C_{0}\omega_{0}^{e}+l)
\end{gather*}
where%
\[
\mathbf{\mu}^{e}=\dfrac{\left(  (1+\rho)m_{1}+\rho m_{2}\right)  }{M_{2}%
}\mathbf{\lambda}^{e},\quad(\omega_{0}^{e})^{2}=\dfrac{(m_{1}+m_{2})h_{1}%
(\rho)}{m_{1}m_{2}}%
\]
is an Eulerian relative equilibrium. The total angular momentum of the system
is expressed by%
\[
\mathbf{L}=(0,0,C_{2}\omega_{2}^{e}+C_{1}\omega_{1}^{e}+C_{0}\omega_{0}%
^{e}+l+g_{1}\omega_{0}^{e}(\lambda^{e})^{2}+g_{2}\omega_{0}^{e}(\mu^{e})^{2})
\]
with $l$ being the gyrostatic momentum. The vectors $\mathbf{\Pi}_{2}^{e},$
$\mathbf{\Pi}_{1}^{e}$ verify the vectorial equations%
\[
\mathbf{\Pi}_{1}^{e}\times\mathbf{\mathbf{\Omega}}_{1}^{e}=\mathbf{0,}%
\quad\mathbf{\Pi}_{2}^{e}\times\mathbf{\mathbf{\Omega}}_{2}^{e}=\mathbf{0}%
\]

\end{proposition}

\subsection{Eulerian relative equilibria when $S_{2}$ and $S_{1}$ are
spherical rigid bodies}

Consider the existence and number of solutions for Eulerian relative
equilibria when $S_{2}$ and $S_{1}$ are spherical rigid bodies. In this case
$C_{1}=A_{1}$, $C_{2}=A_{2}$ and the equation (\ref{EcEuler1}) is equivalent
to the below polynomial equation%
\begin{gather}
(m_{1}+m_{2})\rho^{5}+(3m_{2}+2m_{1})\rho^{4}+(3m_{2}+m_{1}+m_{0}%
(sgn(1+\rho)\medskip\label{quinticaeuler}\\
-sgn(\rho)))\rho^{3}+(m_{2}-m_{2}sgn(1+\rho)-m_{1}sgn(\rho)-3m_{0}%
sgn(\rho))\rho^{2}\medskip\nonumber\\
-(3sgn(\rho)m_{0}+2sgn(\rho)m_{1})\rho-(m_{0}sgn(\rho)+m_{1}sgn(\rho))=0
\end{gather}

This equation has a unique real solution in the intervals $(-\infty,-1)$,
$(-1,0)$ and $(0,+\infty)$. Therefore only one Eulerian relative equilibrium exists.

On the other hand%
\[
\mid\mathbf{\mathbf{\Omega}}_{0}^{e}\mid^{2}=\dfrac{G(m_{1}+m_{2})}%
{\mid\mathbf{\lambda}^{e}\mid^{3}}\left(  1+\dfrac{m_{0}}{(m_{1}+m_{2}%
)}\left(  \dfrac{sgn(1+\rho)}{(1+\rho)^{2}}-\dfrac{sgn(\rho)}{\rho^{2}%
}\right)  \right)
\]
where $\rho$ is the only solution of the equation (\ref{quinticaeuler}).

Proposition 5 gathers the results about Eulerian relative equilibria when
$S_{2}$ and $S_{1}$ are spherical rigid bodies in the configurations $a)$, $b)
$ and $c)$.

\begin{proposition}
\begin{enumerate}
\item If $\rho$ is the unique positive root of the equation%
\begin{gather*}
(m_{1}+m_{2})\rho^{5}+(3m_{2}+2m_{1})\rho^{4}+(3m_{2}+m_{1})\rho^{3}%
+\medskip\\
-(3m_{0}+m_{1})\rho^{2}-(3m_{0}+2m_{1})\rho-(m_{0}+m_{1})=0
\end{gather*}
with
\[
\mid\mathbf{\mathbf{\Omega}}_{0}^{e}\mid^{2}=\dfrac{G(m_{1}+m_{2})}%
{\mid\mathbf{\lambda}^{e}\mid^{3}}\left(  1+\dfrac{m_{0}}{(m_{1}+m_{2}%
)}\left(  \dfrac{1}{(1+\rho)^{2}}-\dfrac{1}{\rho^{2}}\right)  \right)
\]
then $\mathbf{z}_{e}=(\mathbf{\Pi}_{2}^{e},$ $\mathbf{\Pi}_{1}^{e},$
$\mathbf{\Pi}_{0}^{e},$ $\mathbf{\lambda}^{e},$ $\mathbf{p}_{\mathbf{\lambda}%
}^{e},$ $\mathbf{\mu}^{e},$ $\mathbf{p}_{\mu}^{e})$, given by
\begin{gather*}
\mathbf{\lambda}^{e}=(\lambda^{e},0,0),\quad\mathbf{\mu}^{e}=(\mu
^{e},0,0),\quad\mathbf{\mathbf{\Omega}}_{0}^{e}=(0,0,\omega_{0}^{e})\medskip\\
\mathbf{p}_{\mathbf{\lambda}}^{e}=(0,g_{1}\omega_{0}^{e}\lambda^{e}%
,0),\quad\mathbf{p}_{\mathbf{\mu}}^{e}=(0,g_{2}\omega_{0}^{e}\mu^{e}%
,0),\quad\mathbf{\Pi}_{0}^{e}=(0,0,C_{0}\omega_{0}^{e}+l)
\end{gather*}
is the unique solution of relative equilibrium of Euler type in the
configuration $S_{2}S_{1}S_{0}$.

\item If $\rho\in(-1,0)$ is the unique root of the equation%
\begin{gather*}
(m_{1}+m_{2})\rho^{5}+(3m_{2}+2m_{1})\rho^{4}+(3m_{2}+m_{1})\rho^{3}%
+\medskip\\
+(3m_{0}+2m_{2}+m_{1})\rho^{2}+(3m_{0}+2m_{1})\rho+(m_{0}+m_{1})=0
\end{gather*}
with%
\[
\mid\mathbf{\mathbf{\Omega}}_{0}^{e}\mid^{2}=\dfrac{G(m_{1}+m_{2})}%
{\mid\mathbf{\lambda}^{e}\mid^{3}}\left(  1+\dfrac{m_{0}}{(m_{1}+m_{2}%
)}\left(  \dfrac{1}{(1+\rho)^{2}}+\dfrac{1}{\rho^{2}}\right)  \right)
\]
then $\mathbf{z}_{e}=(\mathbf{\Pi}_{2}^{e},$ $\mathbf{\Pi}_{1}^{e},$
$\mathbf{\Pi}_{0}^{e},$ $\mathbf{\lambda}^{e},$ $\mathbf{p}_{\mathbf{\lambda}%
}^{e},$ $\mathbf{\mu}^{e},$ $\mathbf{p}_{\mu}^{e})$ given by
\begin{gather*}
\mathbf{\lambda}^{e}=(\lambda^{e},0,0),\quad\mathbf{\mu}^{e}=(\mu
^{e},0,0),\quad\mathbf{\mathbf{\Omega}}_{0}^{e}=(0,0,\omega_{0}^{e})\medskip\\
\mathbf{p}_{\mathbf{\lambda}}^{e}=(0,g_{1}\omega_{0}^{e}\lambda^{e}%
,0),\quad\mathbf{p}_{\mathbf{\mu}}^{e}=(0,g_{2}\omega_{0}^{e}\mu^{e}%
,0),\quad\mathbf{\Pi}_{0}^{e}=(0,0,C_{0}\omega_{0}^{e}+l)
\end{gather*}
is the unique solution of relative equilibrium of Euler type in the
configuration $S_{2}S_{0}S_{1}.$

\item If $\rho\in(-\infty,-1)\ $is the unique root of the equation
\begin{gather*}
(m_{1}+m_{2})\rho^{5}+(3m_{1}+2m_{2})\rho^{4}+(2m_{0}+3m_{2}+m_{1})\rho
^{3}+\medskip\\
+(3m_{0}+m_{1})\rho^{2}+(3m_{0}+2m_{1})\rho+(m_{0}+m_{1})=0
\end{gather*}
where%
\[
\mid\mathbf{\mathbf{\Omega}}_{0}^{e}\mid^{2}=\dfrac{G(m_{1}+m_{2})}%
{\mid\mathbf{\lambda}^{e}\mid^{3}}\left(  1+\dfrac{m_{0}}{(m_{1}+m_{2}%
)}\left(  \dfrac{1}{\rho^{2}}-\dfrac{1}{(1+\rho)^{2}}\right)  \right)
\]
then $\mathbf{z}_{e}=(\mathbf{\Pi}_{2}^{e},$ $\mathbf{\Pi}_{1}^{e},$
$\mathbf{\Pi}_{0}^{e},$ $\mathbf{\lambda}^{e},$ $\mathbf{p}_{\mathbf{\lambda}%
}^{e},$ $\mathbf{\mu}^{e},$ $\mathbf{p}_{\mu}^{e})$ given by
\begin{gather*}
\mathbf{\lambda}^{e}=(\lambda^{e},0,0),\quad\mathbf{\mu}^{e}=(\mu
^{e},0,0),\quad\mathbf{\mathbf{\Omega}}_{0}^{e}=(0,0,\omega_{0}^{e})\medskip\\
\mathbf{p}_{\mathbf{\lambda}}^{e}=(0,g_{1}\omega_{0}^{e}\lambda^{e}%
,0),\quad\mathbf{p}_{\mathbf{\mu}}^{e}=(0,g_{2}\omega_{0}^{e}\mu^{e}%
,0),\quad\mathbf{\Pi}_{0}^{e}=(0,0,C_{0}\omega_{0}^{e}+l)
\end{gather*}
is the unique solution of relative equilibrium of Euler type in the
configuration $S_{0}S_{2}S_{1}.$
\end{enumerate}
\end{proposition}

\subsection{Eulerian relative equilibria when $S_{2}$ and $S_{1}$ are not
spherical rigid bodies}

In the present case, after carrying out the appropriate calculations, the
equation (\ref{EcEuler1}) is reduced to the study of the positive real roots
of the nine degree equation%
\begin{equation}
\beta_{2}q(\rho)-m_{1}m_{2}a^{2}\rho^{2}(\rho+1)^{2}p(\rho)=0
\label{nueveEuler}%
\end{equation}

where%
\begin{gather}
p(\rho)=(m_{1}+m_{2})\rho^{5}+(3m_{2}+2m_{1})\rho^{4}+(3m_{2}+m_{1}%
+m_{0}(sgn(1+\rho)\medskip\\
-sgn(\rho)))\rho^{3}+(m_{2}-m_{2}sgn(1+\rho)-m_{1}sgn(\rho)-3m_{0}%
sgn(\rho))\rho^{2}\medskip\nonumber\\
-(3sgn(\rho)m_{0}+2sgn(\rho)m_{1})\rho-(m_{0}sgn(\rho)+m_{1}sgn(\rho))
\end{gather}
and%
\begin{gather*}
q(\rho)=m_{0}(ksgn(\rho)m_{2}-sgn(1+\rho)m_{1})\rho^{5}+m_{2}(sgn(1+\rho
)m_{1}+5sgn(\rho)m_{0}k+km_{1}sgn(\rho))\rho^{4}\medskip\\
+2ksgn(\rho)(2m_{1}+5m_{0})\rho^{3}+2ksgn(\rho)(3m_{1}+5m_{0})\rho
^{2}+ksgn(\rho)m_{2}(4m_{1}+5m_{0})\rho\medskip\\
+ksgn(\rho)m_{2}(m_{0}+m_{1})
\end{gather*}
where
\[%
\begin{array}
[c]{cc}%
\beta_{1}=3(C_{1}-A_{1})/2,\ \  & \beta_{2}=3(C_{2}-A_{2})/2
\end{array}
\]
with $\beta_{1}=k\beta_{2},$ $a=\mid\mathbf{\lambda}_{e}\mid$ and
$k\in\mathbb{R}$.

In order to study the number of real roots of the polynomial (\ref{nueveEuler}%
) the rational function will be studied%
\[
\beta_{2}=R(\rho)=\dfrac{m_{1}m_{2}a^{2}\rho^{2}(\rho+1)^{2}p(\rho)}{q(\rho)}%
\]

In practical applications $m_{0}$ is very small, then up to the first order in
$m_{0}$%
\[
\beta_{2}=R_{1}(\rho)=\dfrac{a^{2}\rho^{2}(\rho+1)^{2}p_{1}(\rho)}{q_{1}%
(\rho)}+o(m_{0})
\]
where%
\begin{gather*}
p_{1}(\rho)=\rho^{5}+(2+\mu)\rho^{4}+(1+2\mu)\rho^{3}+(\mu-\mu sgn(1+\rho
)\medskip\\
-(1-\mu)sgn(\rho))\rho^{2}-2sgn(\rho)(1-\mu)\rho-(1-\mu)sgn(\rho)
\end{gather*}
and%
\[
q_{1}(\rho)=(sgn(1+\rho)+sgn(\rho)k)\rho^{4}+4sgn(\rho)k\rho^{3}%
+6sgn(\rho)k\rho^{2}+4sgn(\rho)k\rho+sgn(\rho)k\medskip
\]
where $\mu=\dfrac{m_{1}}{m_{2}+m_{1}}$.

The polynomial $q_{1}$ has no roots in $(0,+\infty)$ and $(-1,0)$ if $k>0$ and
$k<0,$ respectively. On the other hand, $q_{1}$ has only one root, $\rho_{1}$,
in $(0,+\infty)$ and $(-1,0)$ if $k<0$ and $k>0,$ respectively$.$ $\rho_{0}$
will be denoted as the only root of $p_{1}.$ The implicit curve $Res(k,\mu
)=0,$ where $Res$ is the resultant of the polynomials $p_{1}$ and $q_{1}$, is
used to study the graph of $R_{1}.$ When $\mu_{0}$ has a fixed value, the only
$k_{0}$ which verifies $Res(k_{0},$ $\mu_{0})=0,$ exists according to the
Implicit Function Theorem. Consider the only root, $\widetilde{\rho_{1}}%
=\rho_{0}(\mu_{0})$, of $q_{1}$ for $k_{0}$ and $\mu_{0}. $ The expressions%
\[
(\rho_{\max},\text{ }\xi_{1}(k)=R_{1}(\rho_{\max})),\quad(\rho_{\min},\text{
}\xi_{2}(k)=R_{1}(\rho_{\min}))
\]
are the local maximum and minimum of the function $R_{1}.$

In the configuration $S_{2}S_{1}S_{0},$ for any value of $\mu$ fixed%
\[
\lim_{k\rightarrow+\infty}\xi_{2}(k)=0,\quad\lim_{k\rightarrow0^{+}}\xi
_{2}(k)=-\infty
\]
if $k>0.$ For $k<0$%
\[
\lim_{k\rightarrow+\infty}\xi_{1}(k)=0,\quad\lim_{k\rightarrow0^{+}}\xi
_{2}(k)=+\infty
\]
if $\rho_{1}>\widetilde{\rho_{1}}.$ If $\rho_{1}\leq\widetilde{\rho_{1}} $,
$R_{1}$ is strictly increasing.

For the configuration $S_{2}S_{0}S_{1}$ if $\widetilde{\rho_{1}}\leq\rho_{1}$
and $k>0,$ the function $R_{1}$ has just a minimum, which verifies%
\[
\lim_{k\rightarrow0^{+}}\xi_{2}(k)=\xi_{0}%
\]

If $\rho_{1}<\widetilde{\rho_{1}}$, then
\[
\lim_{k\rightarrow+\infty}\xi_{2}(k)=0
\]

If $k<0,$ then $R_{1}$ verifies that%
\begin{align*}
\lim_{k\rightarrow0^{-}}\xi_{2}(k)  &  =\xi_{0},\quad\lim_{k\rightarrow0^{-}%
}\xi_{1}(k)=+\infty\medskip\\
\lim_{k\rightarrow-\infty}\xi_{2}(k)  &  =0,\quad\lim_{k\rightarrow-\infty}%
\xi_{1}(k)=0
\end{align*}
The results for the configuration $S_{0}S_{2}S_{1}$ will be deduced from the
configuration $S_{2}S_{1}S_{0}.$

According to these statements, the following proposition can be stated.

\begin{proposition}
In the configuration $S_{2}S_{1}S_{0}$ the following results are verified.

\begin{description}
\item[a)] For $k>0$ then:
\end{description}

\begin{enumerate}
\item If $\beta_{2}<R_{1}(\rho_{\min})$, then Eulerian relative equilibria do
not exist.

\item If $\beta_{2}=R_{1}(\rho_{\min})$, a unique 2-parametric family of
Eulerian relative equilibria exist.

\item If $R_{1}(\rho_{\min})<\beta_{2}<0,$ two 2-parametric families of
Eulerian relative equilibria exist.

\item If $\beta_{2}>0$ a unique 2-parametric family of Eulerian relative
equilibria exists.
\end{enumerate}

\begin{description}
\item[b)] For $k_{0}<k<0$ and $\beta_{2}>0,$ then:
\end{description}

\begin{enumerate}
\item If $\beta_{2}\in(\xi_{1}(k),$ $\xi_{2}(k)),$ Eulerian relative
equilibria do not exist.

\item If $\beta_{2}=\xi_{1}(k)$ or $\beta_{2}=\xi_{2}(k),$ then a unique
2-parametric family of Eulerian relative equilibria exists.

\item If $\beta_{2}>\xi_{2}(k)$ two 2-parametric families of Eulerian relative
equilibria exists.

\item If $0<\beta_{2}<\xi_{1}(k)$ two 2-parametric families of Eulerian
relative equilibria exists.
\end{enumerate}

\begin{description}
\item[c)] For $k_{0}<k<0$ \ and $\beta_{2}<0,$ then:
\end{description}

\begin{enumerate}
\item A unique 2-parametric family of Eulerian relative equilibria exist.
\end{enumerate}

\begin{description}
\item[d)] For $k<k_{0}$ and $\beta_{2}>0,$ then:
\end{description}

\begin{enumerate}
\item Two 2-parametric families of Eulerian relative equilibria exists.
\end{enumerate}

\begin{description}
\item[e)] For $k<k_{0}$ and $\beta_{2}<0,$ then:
\end{description}

\begin{enumerate}
\item A unique 2-parametric family of Eulerian relative equilibria exist.
\end{enumerate}
\end{proposition}

Similarly we obtain the following result.

\begin{proposition}
In the configuration $S_{2}S_{0}S_{1}$ the following results are verified.

\begin{description}
\item[a)] For $k>0$ then:
\end{description}

\begin{enumerate}
\item If $\beta_{2}<R_{1}(\rho_{\min})$, then a unique 2-parametric family of
Eulerian relative equilibria exists.

\item If $\beta_{2}=R_{1}(\rho_{\min})$, a two 2-parametric family of Eulerian
relative equilibria exist.

\item If $R_{1}(\rho_{\min})<\beta_{2}<0,$ three 2-parametric families of
Eulerian relative equilibria exist.

\item If $\beta_{2}>0$ a unique 2-parametric family of Eulerian relative
equilibria exists.
\end{enumerate}

\begin{description}
\item[b)] For $k<0$ and $\beta_{2}>0,$ then:
\end{description}

\begin{enumerate}
\item If $\beta_{2}\in(\xi_{1}(k),$ $\xi_{2}(k)),$ Eulerian relative
equilibria do not exist.

\item If $\beta_{2}=\xi_{1}(k)$ or $\beta_{2}=\xi_{2}(k)$ a unique
2-parametric family of Eulerian relative equilibria exist.

\item If $\beta_{2}>\xi_{2}(k),$ two 2-parametric families of Eulerian
relative equilibria exist.

\item If $0<\beta_{2}<\xi_{1}(k),$ two 2-parametric families of Eulerian
relative equilibria exist.
\end{enumerate}

\begin{description}
\item[c)] For $k<0$ and $\beta_{2}<0,$ then:
\end{description}

\begin{enumerate}
\item If $R_{1}(\rho_{\min})<\beta_{2}<0,$ then two 2-parametric families of
Eulerian relative equilibria exist.

\item If $R_{1}(\rho_{\min})=\beta_{2},$ then a unique 2-parametric family of
Eulerian relative equilibria exists.

\item If $\beta_{2}<R_{1}(\rho_{\min}),$ then Eulerian relative equilibria do
not exist.

\item If $0<\beta_{2}<R_{2}(\rho_{\max}),$ then two 2-parametric families of
Eulerian relative equilibria exist.

\item If $\beta_{2}=R_{2}(\rho_{\max}),$ then a unique 2-parametric family of
Eulerian relative equilibria exists.

\item If $\beta_{2}>R_{2}(\rho_{\max}),$ then Eulerian relative equilibria do
not exist.
\end{enumerate}
\end{proposition}

\section{Stability of Eulerian relative equilibria}

The tangent flow of the equations (\ref{EcuacHamilt}) in an Eulerian relative
equilibrium $\mathbf{z}_{e}$ is expressed by%
\[
\dfrac{d\delta\mathbf{z}}{dt}=\mathfrak{U(}\mathbf{z}_{e})\delta\mathbf{z}%
\]
with $\delta\mathbf{z=z-z}_{e}$ and $\mathfrak{U(}\mathbf{z}_{e})$ being the
Jacobian matrix of (\ref{EcuacHamilt}) in $\mathbf{z}_{e}.$

The characteristic polynomial $\mathfrak{U(}\mathbf{z}_{e})$ is expressed as
follows%
\begin{equation}
P(\lambda)=\lambda^{3}(\lambda^{2}+\Phi_{0}^{2})(\lambda^{2}+\Phi_{1}%
^{2})(\lambda^{2}+\Phi_{2}^{2})(\lambda^{4}+m\lambda^{2}+n)(\lambda
^{8}+p\lambda^{6}+q\lambda^{4}+r\lambda^{2}+s) \label{polinomiocarac}%
\end{equation}
with $\Phi_{i}^{2}=\dfrac{(C_{i}-A_{i})\omega_{i}^{e}+l}{A_{i}}$. The
coefficients present in the above polynomial are functions of the parameters
of the problem and $\rho$, where $\rho$ is taken as the root of the equation
(\ref{EcEuler1}).

\subsection{$S_{2}$ and $S_{1}$ are spherical rigid bodies}

The characteristic polynomial (\ref{polinomiocarac}) of $\mathfrak{U(}%
\mathbf{z}_{e})$ simplifies to%
\begin{equation}
P(\lambda)=\lambda^{5}(\lambda^{2}+\Phi_{0}^{2})(\lambda^{2}+\Phi_{1}%
^{2})(\lambda^{2}+\Phi_{2}^{2})(\lambda^{2}+(\omega_{0}^{e})^{2})^{2}%
(\lambda^{2}+p)(\lambda^{4}+q\lambda^{2}+r)\nonumber
\end{equation}
with coefficients shown in Appendix B.

If $p\geq0,$ $q\geq0,$ $r\geq0,$ $q^{2}-4r\geq0,$ then $\mathbf{z}_{e}$ is
spectrally stable. These conditions are not verified since $r<0.$

\begin{proposition}
If\textbf{\ }$\mathbf{z}_{e}$ is the only relative equilibrium in the
configuration $S_{0}S_{2}S_{1}$ of the zero order approximate dynamics then it
is unstable.
\end{proposition}

\subsection{$S_{2}$ and $S_{1}$ are close to a sphere}

The case in which $S_{i}$ $(i=1,2)$ are close to a sphere will now be
analyzed. In this case $C_{i}-A_{i}\thickapprox0$, and this is the reason why
by applying the Implicit Function Theorem $\mathbf{z}_{e}$ is unstable.

If $C_{i}-A_{i}$ is not close to zero, the coefficients of the polynomial
(\ref{polinomiocarac}) have very complicated expressions. Numeric calculations
prove that linear stable Eulerian relative equilibria exist for certain values
of the parameters $C_{i}-A_{i}$ $(i=1,2)$ (see Vera and Vigueras \cite{7} for
details). \emph{These results are also applicable to configurations}
$S_{2}S_{0}S_{1}$ \emph{and} $S_{2}S_{1}S_{0}.$

\section{Lagrangian relative equilibria}

\subsection{Necessary conditions of existence}

If $\mathbf{z}_{e}=(\mathbf{\Pi}_{2}^{e},$ $\mathbf{\Pi}_{1}^{e},$
$\mathbf{\Pi}_{0}^{e},$ $\mathbf{\lambda}^{e},$ $\mathbf{p}_{\mathbf{\lambda}%
}^{e},$ $\mathbf{\mu}^{e},$ $\mathbf{p}_{\mu}^{e})$ is a Lagrangian relative
equilibrium, the below identities are verified%
\begin{gather*}
\mathbf{\lambda}^{e}\ \mathbf{\times\ }(\mathbf{\nabla}_{\mathbf{\lambda}%
}\mathcal{V})_{e}=\mathbf{0,}\quad g_{1}\mid\mathbf{\mathbf{\Omega}}_{0}%
^{e}\mid^{2}(\mathbf{\lambda}^{e}\ \mathbf{\times\ \mu}^{e})=(\mathbf{\nabla
}_{\mathbf{\lambda}}\mathcal{V})_{e}\ \mathbf{\times\ \mu}^{e}\medskip\\
\mathbf{\mu}^{e}\ \mathbf{\times\ }(\mathbf{\nabla}_{\mathbf{\mu}}%
\mathcal{V})_{e}=\mathbf{0,}\quad g_{2}\mid\mathbf{\mathbf{\Omega}}_{0}%
^{e}\mid^{2}(\mathbf{\lambda}^{e}\ \mathbf{\times\ \mu}^{e})=\mathbf{\lambda
}^{e}\ \mathbf{\times\ }(\mathbf{\nabla}_{\mathbf{\mu}}\mathcal{V})_{e}%
\end{gather*}

After some calculations, it is concluded%
\begin{gather*}
(A_{12})_{e}(\mathbf{\lambda}^{e}\ \mathbf{\times\ \mu}^{e})=\mathbf{0,}\quad
g_{1}\mid\mathbf{\mathbf{\Omega}}_{0}^{e}\mid^{2}(\mathbf{\lambda}%
^{e}\ \mathbf{\times\ \mu}^{e})=(A_{11})_{e}(\mathbf{\lambda}^{e}%
\ \mathbf{\times\ \mu}^{e})\medskip\\
(A_{21})_{e}(\mathbf{\lambda}^{e}\ \mathbf{\times\ \mu}^{e})=\mathbf{0,}\quad
g_{2}\mid\mathbf{\mathbf{\Omega}}_{0}^{e}\mid^{2}(\mathbf{\lambda}%
^{e}\ \mathbf{\times\ \mu}^{e})=(A_{22})_{e}(\mathbf{\lambda}^{e}%
\ \mathbf{\times\ \mu}^{e})
\end{gather*}

The expressions $(A_{ij})_{e}$ are certain values of the potential function in
the Lagrangian relative equilibrium (see Appendix A, for details).

Concluding, the following relationships are verified%
\[
(A_{12})_{e}=0,\quad\mid\mathbf{\mathbf{\Omega}}_{0}^{e}\mid^{2}%
=\dfrac{(A_{11})_{e}}{g_{1}}=\dfrac{(A_{22})_{e}}{g_{2}}%
\]
and the below proposition can be obtained.

\begin{proposition}
Let $\mathbf{z}_{e}=(\mathbf{\Pi}_{2}^{e},$ $\mathbf{\Pi}_{1}^{e},$
$\mathbf{\Pi}_{0}^{e},$ $\mathbf{\lambda}^{e},$ $\mathbf{p}_{\mathbf{\lambda}%
}^{e},$ $\mathbf{\mu}^{e},$ $\mathbf{p}_{\mu}^{e})$ be a Lagrangian relative
equilibrium. Then it is seen that%
\begin{align*}
g_{2}(\widetilde{A}_{11})_{e}  &  =g_{1}(\widetilde{A}_{22})_{e}%
,\quad(\widetilde{A}_{12})_{e}=0,\quad\mid\mathbf{\mathbf{\Omega}}_{0}^{e}%
\mid^{2}=\dfrac{(A_{11})_{e}}{g_{1}}\\
&
\end{align*}

\end{proposition}

If $\mathbf{z}_{e}=(\mathbf{\Pi}_{2}^{e},$ $\mathbf{\Pi}_{1}^{e},$
$\mathbf{\Pi}_{0}^{e},$ $\mathbf{\lambda}^{e},$ $\mathbf{p}_{\mathbf{\lambda}%
}^{e},$ $\mathbf{\mu}^{e},$ $\mathbf{p}_{\mu}^{e})$ is a Lagrangian relative
equilibrium, then using the expressions of $(A_{ij})_{e}$ it is proven%
\[
\dfrac{m_{1}}{\mid\mathbf{\lambda}^{e}\mid^{3}}=\dfrac{m_{1}}{\mid\mathbf{\mu
}^{e}\mathbf{-}\frac{m_{2}}{M_{2}}\mathbf{\lambda}^{e}\mid^{3}}+\dfrac
{\beta_{1}}{\mid\mathbf{\mu}^{e}\mathbf{-}\frac{m_{2}}{M_{2}}\mathbf{\lambda
}^{e}\mid^{5}},\text{ \ }\dfrac{m_{2}}{\mid\mathbf{\lambda}^{e}\mid^{3}%
}=\dfrac{m_{2}}{\mid\mathbf{\mu}^{e}\mathbf{+}\frac{m_{1}}{M_{2}%
}\mathbf{\lambda}^{e}\mid^{3}}+\dfrac{\beta_{2}}{\mid\mathbf{\mu}%
^{e}\mathbf{+}\frac{m_{1}}{M_{2}}\mathbf{\lambda}^{e}\mid^{5}}%
\]
and
\[
\mid\mathbf{\mathbf{\Omega}}_{0}^{e}\mid^{2}=\dfrac{G(m_{0}+m_{1}+m_{2}%
)}{Z^{3}}%
\]

From the previous expressions, the bellow result is deduced.

\begin{proposition}
If $\mathbf{z}_{e}=(\mathbf{\Pi}_{2}^{e},$ $\mathbf{\Pi}_{1}^{e},$
$\mathbf{\Pi}_{0}^{e},$ $\mathbf{\lambda}^{e},$ $\mathbf{p}_{\mathbf{\lambda}%
}^{e},$ $\mathbf{\mu}^{e},$ $\mathbf{p}_{\mu}^{e})$ is a Lagrangian relative
equilibrium, then by denoting $\mid\mathbf{\lambda}^{e}\mid=Z,$ $\mid
\mathbf{\mu}^{e}\mathbf{-}\frac{m_{2}}{M_{2}}\mathbf{\lambda}^{e}\mid=X,$
$\mid\mathbf{\mu}^{e}\mathbf{+}\frac{m_{1}}{M_{2}}\mathbf{\lambda}^{e}\mid=Y,$
the system of equations
\begin{equation}
\left\{
\begin{array}
[c]{c}%
X^{5}-Z^{3}X^{2}-\beta_{1}Z^{3}=0\medskip\\
Y^{5}-Z^{3}Y^{2}-\beta_{2}Z^{3}=0
\end{array}
\right.  \label{Lagrange2}%
\end{equation}
has positive real solutions.
\end{proposition}

\begin{remark}
$Z$, $\beta_{1}$ and $\beta_{2}$ will be the parameters that exert a strong
influence on the study of the number of the different configurations of a
Lagrangian relative equilibrium.
\end{remark}

\subsection{Sufficient conditions of existence}

If $Z$ has a fixed value and $X$, $Y$ verify the below system of equations%
\begin{equation}
\left\{
\begin{array}
[c]{c}%
X^{5}-Z^{3}X^{2}-\beta_{1}Z^{3}=0\medskip\\
Y^{5}-Z^{3}Y^{2}-\beta_{2}Z^{3}=0
\end{array}
\right.
\end{equation}
, with respect to an appropriate reference system, the coordinates of
Lagrangian relative equilibria are completely determinated. If $X=Y\neq Z$ is
a solution to the previous system, then the gyrostat $S_{0}$ and the rigid
bodies $S_{i}$ $(i=1,2)$ form an isosceles triangle when $\beta_{1}=\beta_{2}%
$. On the other hand, if $\beta_{1}\neq\beta_{2}$ then $X\neq Y\neq Z$ and
$S_{i}$ $(i=0,1,2)$ form a scalene triangle.

The following proposition, whose demonstration is based on verifying the
equations of the equilibria, shows how Lagrangian relative equilibria are when
$S_{0}$, $S_{1},$ $S_{2}$ form an isosceles triangle. Following the same line
of reasoning, Lagrangian relative equilibria could be described when $S_{0}$,
$S_{1},$ $S_{2}$ form a scalene triangle.

\begin{proposition}
With respect to an appropriate reference system, $\mathbf{z}_{e}=(\mathbf{\Pi
}_{2}^{e},$ $\mathbf{\Pi}_{1}^{e},$ $\mathbf{\Pi}_{0}^{e},$ $\mathbf{\lambda
}^{e},$ $\mathbf{p}_{\mathbf{\lambda}}^{e},$ $\mathbf{\mu}^{e},$
$\mathbf{p}_{\mu}^{e})$ expressed by%
\begin{gather}
\mathbf{\lambda}^{e}=(x_{1},y_{1},0),\quad\mathbf{p}_{\mathbf{\lambda}}%
^{e}=g_{1}\omega_{0}^{e}(-y_{1},x_{1},0),\quad\mathbf{\mu}^{e}=(x_{2}%
,y_{2},0)\medskip\label{parametri}\\
\mathbf{p}_{\mathbf{\mu}}^{e}=g_{2}\omega_{0}^{e}(-y_{2},x_{2},0),\quad
\mathbf{\mathbf{\Omega}}_{0}^{e}=(0,0,\omega_{0}^{e}),\quad\mathbf{\Pi}%
_{0}^{e}=(0,0,C_{0}\omega_{0}^{e}+l)\nonumber
\end{gather}
with
\[
x_{1}=Z,\text{ \ }y_{1}=0,\text{ \ }x_{2}=\dfrac{Z(m_{2}-m_{1})}{2(m_{1}%
+m_{2})},\text{ \ }y_{2}=\pm\dfrac{\sqrt{4X^{2}-Z^{2}}}{2},\text{ \ }%
(\omega_{0}^{e})^{2}=\dfrac{G(m_{0}+m_{1}+m_{2})}{Z^{3}}%
\]
are an isosceles Lagrangian relative equilibrium. The total angular momentum
vector of the system is given by%
\[
\mathbf{L}=(0,0,C_{2}\omega_{2}^{e}+C_{1}\omega_{1}^{e}+C_{0}\omega_{0}%
^{e}+l+\omega_{0}^{e}%
{\displaystyle\sum\limits_{i=1}^{2}}
g_{i}(x_{i}^{2}+y_{i}^{2}))
\]

\end{proposition}

\begin{proposition}
With respect to an appropriate reference system, $\mathbf{z}_{e}=(\mathbf{\Pi
}_{2}^{e},$ $\mathbf{\Pi}_{1}^{e},$ $\mathbf{\Pi}_{0}^{e},$ $\mathbf{\lambda
}^{e},$ $\mathbf{p}_{\mathbf{\lambda}}^{e},$ $\mathbf{\mu}^{e},$
$\mathbf{p}_{\mu}^{e})$ given by%
\begin{gather}
\mathbf{\lambda}^{e}=(x_{1},y_{1},0),\quad\mathbf{p}_{\mathbf{\lambda}}%
^{e}=g_{1}\omega_{0}^{e}(-y_{1},x_{1},0),\quad\mathbf{\mu}^{e}=(x_{2}%
,y_{2},0)\medskip\\
\mathbf{p}_{\mathbf{\mu}}^{e}=g_{2}\omega_{0}^{e}(-y_{2},x_{2},0),\quad
\mathbf{\mathbf{\Omega}}_{0}^{e}=(0,0,\omega_{0}^{e}),\quad\mathbf{\Pi}%
_{0}^{e}=(0,0,C_{0}\omega_{0}^{e}+l)\nonumber
\end{gather}
with
\begin{gather*}
x_{1}=Z,\text{ \ }y_{1}=0,\text{ \ }x_{2}=\dfrac{m_{2}(X^{2}+Z^{2}%
-Y^{2})-m_{1}(Y^{2}+Z^{2}-X^{2})}{2(m_{1}+m_{2})Z},\text{ }\\
\text{\ }(\omega_{0}^{e})^{2}=\dfrac{G(m_{0}+m_{1}+m_{2})}{Z^{3}},\text{
}y_{2}=\pm\dfrac{\sqrt{(Z+X+Y)(Z+X-Y)(Z+Y-X)(X+Y-Z)}}{2Z}%
\end{gather*}
are a scalene Lagrangian relative equilibrium. The total angular momentum
vector of the system is expressed as follows%
\[
\mathbf{L}=(0,0,C_{2}\omega_{2}^{e}+C_{1}\omega_{1}^{e}+C_{0}\omega_{0}%
^{e}+l+\omega_{0}^{e}%
{\displaystyle\sum\limits_{i=1}^{2}}
g_{i}(x_{i}^{2}+y_{i}^{2}))
\]

\end{proposition}

Next the Lagrangian relative equilibria will be studied when $S_{2}$ and
$S_{1}$ are spherical rigid bodies.

\subsection{Lagrange relative equilibria when $S_{2}$ and $S_{1}$ are
spherical rigid bodies}

If $S_{2}$ and $S_{1}$ are spherical rigid bodies,$\ $then $C_{1}=A_{1}$,
$C_{2}=A_{2}.$ The equations (\ref{Lagrange2}) are%
\[
\left\{
\begin{array}
[c]{c}%
X^{3}=Z^{3}\medskip\\
Y^{3}=Z^{3}%
\end{array}
\right.
\]
being easily deduced from them that $X=Y=Z$, i.e. $S_{0}$, $S_{1}$ and $S_{2}
$ form an equilateral triangle. It is also obtained%
\[
\mid\mathbf{\mathbf{\Omega}}_{0}^{e}\mid^{2}=\dfrac{G(m_{0}+m_{1}+m_{2}%
)}{Z^{3}}%
\]

On the other hand, one parametrization of the former equilibria $\mathbf{z}%
_{e}=(\mathbf{\Pi}_{2}^{e},$ $\mathbf{\Pi}_{1}^{e},$ $\mathbf{\Pi}_{0}^{e},$
$\mathbf{\lambda}^{e},$ $\mathbf{p}_{\mathbf{\lambda}}^{e},$ $\mathbf{\mu}%
^{e},$ $\mathbf{p}_{\mu}^{e})$ is given by the relationship (\ref{parametri})
being%
\[
x_{1}=Z,\text{ \ }y_{1}=0,\text{ \ }x_{2}=\dfrac{Z(m_{2}-m_{1})}{2(m_{1}%
+m_{2})},\text{ \ }y_{2}=y_{2}=\pm\dfrac{\sqrt{3}Z}{2},\text{ \ }(\omega
_{0}^{e})^{2}=\dfrac{G(m_{0}+m_{1}+m_{2})}{Z^{3}}%
\]

This parametrization of the relative equilibria will be useful to study the
stability of the same ones.

\subsection{Lagrangian relative equilibria when $S_{2}$ and $S_{1}$ are not
spherical rigid bodies}

The number of real positive solutions of the below system will be now studied%
\begin{equation}
\left\{
\begin{array}
[c]{c}%
X^{5}-Z^{3}X^{2}-\beta_{1}Z^{3}=0\medskip\\
Y^{5}-Z^{3}Y^{2}-\beta_{2}Z^{3}=0
\end{array}
\right.  \label{EAlgLagrange1}%
\end{equation}
where $Z$, $\beta_{1}$ and $\beta_{2}$ are parameters. With this aim in view,
let us proceed with the study of the number of the positive real roots of the
polynomial%
\[
p(X)=X^{5}-Z^{3}X^{2}-\beta Z^{3}%
\]
according to the values of the parameters $Z$ and $\beta$.

If $\beta\geq0,$ then this polynomial can only have a positive real root by
applying the rule of the signs of Descartes.

If $\beta<0,$ then we can have two positive real roots, just one positive real
root or none. The discriminant of the polynomial, denoted by $discrim(p,X)$,
is expressed by%
\[
discrim(p,X)=\beta Z^{12}(3125\beta^{3}+108Z^{6})
\]

Then if $discrim(p,X)<0$ the polynomial $p$ has two real roots, if
$discrim(p,x)=0$ it has a positive double root, whereas if $discrim(p,x)>0$ it
has no positive root.

The discriminant is zero when the following relationship is verified%
\[
\beta=-\dfrac{3\sqrt[3]{20}}{25}Z^{2}%
\]

From the previous results a detailed study can be carried out of the
bifurcations of the Lagrangian relative equilibria when $\beta_{1}=\beta_{2}.
$

\begin{proposition}
Let $\mathbf{z}_{e}=(\mathbf{\Pi}_{2}^{e},$ $\mathbf{\Pi}_{1}^{e},$
$\mathbf{\Pi}_{0}^{e},$ $\mathbf{\lambda}^{e},$ $\mathbf{p}_{\mathbf{\lambda}%
}^{e},$ $\mathbf{\mu}^{e},$ $\mathbf{p}_{\mu}^{e})$ be a Lagrangian relative
equilibrium and $\beta_{1}=\beta_{2}$, then

\begin{enumerate}
\item If $\beta_{1}>0$, a unique 2-parametric family exists, making $S_{0}$,
$S_{1}$ and $S_{2}$ an isosceles triangle.

\item If $\beta_{1}<0$ (gyrostat prolate), then:

\begin{description}
\item[a1)] If $\dfrac{-7Z^{2}}{32}<\beta_{1}<0$, there are two types of
relative equilibria:

\begin{itemize}
\item One 2-parametric family of relative equilibria, making $S_{0}$, $S_{1}$
and $S_{2}$ an isosceles triangle with $X=Y\neq Z.$

\item Two 2-parametric families of relative equilibria, making $S_{0}$,
$S_{1}$ and $S_{2}$ a scalene triangle with $X\neq Y\neq Z.$
\end{itemize}

\item[a2)] If $-\dfrac{3\sqrt[3]{20}}{25}Z^{2}<\beta_{1}<\dfrac{-7Z^{2}}{32}
$, there are two types of relative equilibria:

\begin{itemize}
\item Two 2-parametric families of relative equilibria, making $S_{0}$,
$S_{1}$ and $S_{2}$ an isosceles triangle with $X=Y\neq Z.$

\item Four 2-parametric families of relative equilibria, making $S_{0}$,
$S_{1}$ and $S_{2}$ an scalene triangle with $X\neq Y\neq Z.$
\end{itemize}

\item[b)] If $\beta_{1}=-\dfrac{3\sqrt[3]{20}}{25}Z^{2},$ a unique
2-parametric family exists making $S_{0}$, $S_{1}$ and $S_{2}$ an isosceles
triangle, with $X=Y\neq Z.$

\item[c)] If $\beta_{1}<-\dfrac{3\sqrt[3]{20}}{25}Z^{2}$, relative equilibria
do not exist.
\end{description}
\end{enumerate}
\end{proposition}

On the other hand, we have that $\mathbf{z}_{e}=(\mathbf{\Pi}_{2}^{e},$
$\mathbf{\Pi}_{1}^{e},$ $\mathbf{\Pi}_{0}^{e},$ $\mathbf{\lambda}^{e},$
$\mathbf{p}_{\mathbf{\lambda}}^{e},$ $\mathbf{\mu}^{e},$ $\mathbf{p}_{\mu}%
^{e}),$ with equations (\ref{parametri}) given by%
\[
x_{1}=Z,\text{ \ }y_{1}=0,\text{ \ }x_{2}=\dfrac{Z(m_{2}-m_{1})}{2(m_{1}%
+m_{2})},\text{ \ }y_{2}=\pm\dfrac{\sqrt{4X^{2}-Z^{2}}}{2},\text{ \ }%
(\omega_{0}^{e})^{2}=\dfrac{G(m_{0}+m_{1}+m_{2})}{Z^{3}}%
\]
describes a parametrization of the Lagrangian relative equilibria when $S_{0}
$, $S_{1}$ and $S_{2}$ form an isosceles triangle. Similar results are
obtained when $S_{0}$, $S_{1}$ and $S_{2}$ form a scalene triangle.

\begin{proposition}
Let $\mathbf{z}_{e}=(\mathbf{\Pi}_{2}^{e},$ $\mathbf{\Pi}_{1}^{e},$
$\mathbf{\Pi}_{0}^{e},$ $\mathbf{\lambda}^{e},$ $\mathbf{p}_{\mathbf{\lambda}%
}^{e},$ $\mathbf{\mu}^{e},$ $\mathbf{p}_{\mu}^{e})$ be a Lagrangian relative
equilibrium and $\beta_{1}\neq\beta_{2}$, then

\begin{enumerate}
\item If $\beta_{1},$ $\beta_{2}>0$ a unique 3-parametric family exists,
making $S_{0}$, $S_{1}$ and $S_{2}$ a scalene triangle.

\item If $\beta_{1}<0$ and $\beta_{2}>0$ then

\begin{description}
\item[a1)] If $\dfrac{-7Z^{2}}{32}<\beta_{1}<0$, there are two 3-parametric
families of relative equilibria, making $S_{0}$, $S_{1}$ and $S_{2}$ a scalene
triangle with $X\neq Y\neq Z.$

\item[a2)] If $-\dfrac{3\sqrt[3]{20}}{25}Z^{2}<\beta_{1}<\dfrac{-7Z^{2}}{32}
$, there are four 3-parametric families of relative equilibria, making $S_{0}
$, $S_{1}$ and $S_{2}$ a scalene triangle with $X\neq Y\neq Z.$

\item[b)] If $\beta_{1}=-\dfrac{3\sqrt[3]{20}}{25}Z^{2}$ a unique 3-parametric
family exists, making $S_{0}$, $S_{1}$ and $S_{2}$ a scalene triangle with
$X=Y\neq Z.$

\item[c)] If $\beta_{1}<-\dfrac{3\sqrt[3]{20}}{25}Z^{2},$ relative equilibria
do not exist.
\end{description}

\item The case $\beta_{2}<0$ and $\beta_{1}>0$ is similar to the previous one.

\item If $\beta_{1},$ $\beta_{2}<0$ then:

\begin{description}
\item[b1)] If $\dfrac{-7Z^{2}}{32}<\beta_{1}<0$ and $\dfrac{-7Z^{2}}{32}%
<\beta_{2}<0,$ there are four 3-parametric families of relative equilibria,
making $S_{0}$, $S_{1}$ and $S_{2}$ a scalene triangle with $X\neq Y\neq Z.$

\item[b2)] If $\dfrac{-7Z^{2}}{32}<\beta_{1}<0$ and $-\dfrac{3\sqrt[3]{20}%
}{25}Z^{2}<\beta_{2}<\dfrac{-7Z^{2}}{32},$ there are six 3-parametric families
of relative equilibria, making $S_{0}$, $S_{1}$ and $S_{2}$ a scalene triangle
with $X\neq Y\neq Z.$

\item[b3)] If $-\dfrac{3\sqrt[3]{20}}{25}Z^{2}<\beta_{1}<\dfrac{-7Z^{2}}{32} $
and $-\dfrac{3\sqrt[3]{20}}{25}Z^{2}<\beta_{2}<\dfrac{-7Z^{2}}{32},$ there are
eight 3-parametric families of relative equilibria, making $S_{0}$, $S_{1}$
and $S_{2}$ a scalene triangle with $X\neq Y\neq Z.$

\item[b4)] If $\beta_{1}=-\dfrac{3\sqrt[3]{20}}{25}Z^{2}$ and $-\dfrac
{3\sqrt[3]{20}}{25}Z^{2}<\beta_{2}<\dfrac{-7Z^{2}}{32},$ there are two
3-parametric families of relative equilibria, making $S_{0}$, $S_{1}$ and
$S_{2}$ a scalene triangle with $X\neq Y\neq Z.$

\item[b5)] If $\beta_{1}=-\dfrac{3\sqrt[3]{20}}{25}Z^{2}$ and $\beta
_{2}<-\dfrac{3\sqrt[3]{20}}{25}Z^{2}$ or $\beta_{2}=-\dfrac{3\sqrt[3]{20}}%
{25}Z^{2} $ and $\beta_{1}<-\dfrac{3\sqrt[3]{20}}{25}Z^{2},$ relative
equilibria do not exist.
\end{description}
\end{enumerate}
\end{proposition}


\subsection{Lagrangian relative equilibria when $S_{2}$ and $S_{1}$ are close
to spherical rigid bodies}

If $\beta_{1},$ $\beta_{2}\approx0,$ the solutions of (\ref{EAlgLagrange1}) up
to the second order in the parameters $\beta_{1}$ and $\beta_{2}$ are%
\begin{gather*}
X=Z+\dfrac{\beta_{1}}{3Z}-\dfrac{\beta_{1}^{2}}{3Z^{3}}+o(\beta_{1}^{2})\\
Y=Z+\dfrac{\beta_{2}}{3Z}-\dfrac{\beta_{2}^{2}}{3Z^{3}}+o(\beta_{2}^{2})
\end{gather*}

Using the previous relationships, the expressions of $\mathbf{\mu}^{e}%
=(x_{2},y_{2},0)$ are obtained with%
\begin{gather*}
x_{2}=\dfrac{(m_{2}-m_{1})Z}{2(m_{2}+m_{1})}-\dfrac{\beta_{1}}{3Z}%
+\dfrac{\beta_{2}}{3Z}+\dfrac{5\beta_{1}^{2}}{18Z^{3}}-\dfrac{5\beta_{2}^{2}%
}{18Z^{3}}+o(\beta_{1}^{2}+\beta_{2}^{2})\medskip\\
y_{2}=\dfrac{\sqrt{3}Z}{2}+\dfrac{\sqrt{3}}{9}(\beta_{1}+\beta_{2}%
)-(\dfrac{23\sqrt{3}\beta_{1}^{2}}{162Z^{3}}-\dfrac{4\sqrt{3}\beta_{1}%
\beta_{2}}{81Z^{3}}+\dfrac{23\sqrt{3}\beta_{2}^{2}}{162Z^{3}})+o(\beta_{1}%
^{2}+\beta_{2}^{2})
\end{gather*}
and $\mathbf{\lambda}^{e}=(Z,0,0)$.

\section{Linear stability of the Lagrangian relative equilibria}

In the equilibrium $\mathbf{z}_{e}$ the tangent flow of the equations
(\ref{EcuacHamilt}) are%
\[
\dfrac{d\delta\mathbf{z}}{dt}=\mathfrak{U(}\mathbf{z}_{e})\delta\mathbf{z}%
\]
where $\delta\mathbf{z=z-z}_{e}$ and $\mathfrak{U(}\mathbf{z}_{e})$ is the
Jacobian matrix of (\ref{EcuacHamilt}) in $\mathbf{z}_{e}.$

\subsection{$S_{2}$ and $S_{1}$ are spherical rigid bodies}

The characteristic polynomial of $\mathfrak{U(}\mathbf{z}_{e})$ has the bellow
expression%
\[
P(\lambda)=\lambda^{5}(\lambda^{2}+\Phi_{0}^{2})(\lambda^{2}+\Phi_{1}%
^{2})(\lambda^{2}+\Phi_{2}^{2})(\lambda^{2}+\omega_{e}^{2})^{3}(\lambda
^{4}+\omega_{e}^{2}\lambda^{2}+q)
\]
where%
\[
\omega_{e}^{2}=\dfrac{G(m_{0}+m_{1}+m_{2})}{Z^{3}},\quad q=\dfrac
{27G^{2}(m_{1}m_{0}+m_{2}m_{0}+m_{1}m_{2})}{4Z^{6}}%
\]
and $\Phi_{i}^{2}=\dfrac{(C_{i}-A_{i})\omega_{i}^{e}+l}{A_{i}}.$

The minimum polynomial of $\mathfrak{U(}\mathbf{z}_{e})$ has this expression%
\[
Q(\lambda)=\lambda^{2}(\lambda^{2}+\Phi_{0}^{2})(\lambda^{2}+\Phi_{1}%
^{2})(\lambda^{2}+\Phi_{2}^{2})(\lambda^{2}+\omega_{e}^{2})(\lambda^{4}%
+\omega_{e}^{2}\lambda^{2}+q)
\]

Then the bellow results are verified.

\begin{proposition}
$\mathbf{z}_{e}$ is spectral stable if%
\begin{equation}
(m_{0}+m_{2}+m_{1})^{2}\geq27(m_{1}m_{0}+m_{2}m_{0}+m_{1}m_{2})
\label{stability1}%
\end{equation}
If%
\[
(m_{0}+m_{2}+m_{1})^{2}<27(m_{1}m_{0}+m_{2}m_{0}+m_{1}m_{2})
\]
then $\mathbf{z}_{e}$ is unstable
\end{proposition}

As the minimum polynomial of $\mathfrak{U(}\mathbf{z}_{e})$ has the
$\lambda=0$ as double root, the matrix $\mathfrak{U(}\mathbf{z}_{e})$ is not
diagonalizable, being the following proposition verified.

\begin{proposition}
The linear system%
\[
\dfrac{d\delta\mathbf{z}}{dt}=\mathfrak{U(}\mathbf{z}_{e})\delta\mathbf{z}%
\]
is unstable.
\end{proposition}

\subsection{$S_{2}$ and $S_{1}$ are not spherical rigid bodies}

Similar results show that the characteristic polynomial has this expression%
\[
P(\lambda)=\lambda^{3}(\lambda^{2}+\Phi_{0}^{2})(\lambda^{2}+\Phi_{1}%
^{2})(\lambda^{2}+\Phi_{2}^{2})(\lambda^{2}+m)(\lambda^{2}+n)(\lambda
^{8}+p\lambda^{6}+q\lambda^{4}+r\lambda^{2}+s)
\]

The minimum polynomial of $\mathfrak{U(}\mathbf{z}_{e})$ is expressed by%
\[
Q(\lambda)=\lambda(\lambda^{2}+\Phi_{0}^{2})(\lambda^{2}+\Phi_{1}^{2}%
)(\lambda^{2}+\Phi_{2}^{2})(\lambda^{2}+m)(\lambda^{2}+n)(\lambda^{8}%
+p\lambda^{6}+q\lambda^{4}+r\lambda^{2}+s)
\]

The coefficients of the characteristic polynomial will be expressed in
function of the parameters of the problem, i.e. the masses and the
coefficients $\beta_{i}$ $(i=1,2)$. We conclude the following result by
employing the Sturm Theorem.

We conclude that the Lagrangian relative equilibria are spectral stable
(lineally stable) if the following conditions are verified%
\begin{gather*}
p^{2}q^{2}-3rp^{3}-6p^{2}s-4q^{3}+14pqr+16qs-18r^{2}\geq0\text{ }%
(>0)\medskip\\
p^{2}qr-48sr-9sp^{3}+32pqs-4q^{2}r+3pr^{2}\geq0\text{ }(>0)\medskip\\
r,s\geq0\text{ }(>0),\text{ }3p^{2}-8q\geq0\text{ }(>0),\text{ }%
pr-16s\geq0\text{ }(>0)\medskip\\
m,n\geq0\text{ }(>0),\text{ \ \ }discrim(h)\geq0\text{ }(>0)\medskip
\end{gather*}

If $S_{i}$ $(i=1,2)$ are arbitrary rigid bodies, the former conditions are
very complicated expressions in the parameters of the problem and can only be
studied by means of numerical analysis.

If $S_{i}$ $(i=1,2)$ are close to a sphere, the coefficients of $P$, up to the
first order in the parameters $\beta_{1}$ and $\beta_{2}$ are expressed as
follows%
\begin{align*}
m  &  =\dfrac{G(m_{0}+m_{1}+m_{2})}{Z^{3}}+o(\beta_{1})+o(\beta_{2}),\quad
n=\dfrac{G(m_{0}+m_{1}+m_{2})}{Z^{3}}+o(\beta_{1})+o(\beta_{2})\medskip\\
r  &  =\dfrac{27G^{3}(m_{0}+m_{1}+m_{2})(m_{1}m_{0}+m_{2}m_{0}+m_{1}m_{2}%
)}{4Z^{9}}+o(\beta_{1})+o(\beta_{2})\medskip\\
q  &  =\dfrac{G^{2}(4m_{0}^{2}+4m_{1}^{2}+4m_{2}^{2}+35m_{0}m_{1}+35m_{0}%
m_{2}+35m_{1}m_{2})}{4Z^{6}}+o(\beta_{1})+o(\beta_{2})\medskip\\
p  &  =\dfrac{2G(m_{0}+m_{1}+m_{2})}{Z^{3}}+o(\beta_{1})+o(\beta_{2}),\quad
s=o(\beta_{1})+o(\beta_{2})
\end{align*}
If the function
\[
s=\dfrac{81G^{4}m_{0}(m_{0}+m_{1}+m_{2})^{2}(\beta_{1}m_{1}+\beta_{2}m)}%
{4}+o(\beta_{1}^{2}+\beta_{2}^{2})
\]
is positive and%
\begin{equation}
(m_{0}+m_{2}+m_{1})^{2}>27(m_{1}m_{0}+m_{2}m_{0}+m_{1}m_{2})\nonumber
\end{equation}
$\mathbf{z}_{e}$ is linearly stable. Then if $S_{i}$ $(i=1,2)$ is close to a
sphere and%
\[
m_{1}(C_{1}-A_{1})+m_{2}(C_{2}-A_{2})>0
\]
$\mathbf{z}_{e}$ is linearly stable when (\ref{stability1}) is verified. If
$m_{1}(C_{1}-A_{1})+m_{2}(C_{2}-A_{2})=0,$ second order terms of the
coefficients in $\beta_{1}$ and $\beta_{2}$ are necessary for the study of the problem.

\section{Conclusions and future works}

In this paper it has been investigated some important periodic solutions of
the dynamics of a gyrostat in Newtonian interaction with two symmetric rigid
bodies. With the hypotheses formulated in the introduction of this paper,
working in the double reduced space of configuration of the problem, both the
equations of motion and the those which determine the relative equilibria have
been derived.

Two families of relative equilibria, \emph{Eulerian} and \emph{Lagrangian} are
studied. The Eulerian relative equilibria have been completely determinated by
a polynomial equation of degree nine. The bifurcations of these equilibria
have been carried out when $m_{0}$ is very small. The bifurcations of the
Lagrangian relative equilibria have been completely investigated and
expressions for the Lagrangian relative equilibria, when both solids are close
to a sphere, are provided. These expressions are useful for the posterior
study of the stability of the same ones.

The instability of Eulerian relative has been proven and necessary and
sufficient conditions for lineal stability of Lagrangian relative equilibria
are provided.

Different results, which had been obtained by means of classic methods in
previous works, have been generalized in a different way, and other results,
not previously considered, have been studied. Some interesting numeric results
have been detailed in Appendix B.

Numerous problems are open, and among them the study of the "inclined"
relative equilibria, in which $\mathbf{\Omega}_{0}^{e}$ form an angle
$\alpha\neq0,$ $\pi/2$ with the vector $\mathbf{\lambda}^{e}\mathbf{\times}$
$\mathbf{\mu}^{e},$ is considered.

The methods used in this work are susceptible of being used in similar
problems. The nonlinear stability of the relative equilibria obtained here is
the logical continuation of this work.

\section{The $(A_{ij})_{e},$ $i,j=1,2.$}

The bellow expressions of the potential $\mathcal{V}$ are obtained.%
\begin{gather*}
(\mathbf{\nabla}_{\mathbf{\lambda}}\mathcal{V})_{e}=\dfrac{Gm_{1}%
m_{2}\mathbf{\lambda}^{e}}{\mid\mathbf{\lambda}^{e}\mid^{3}}-\dfrac
{Gm_{0}m_{2}}{M_{2}}%
{\displaystyle\sum\limits_{i=0}^{1}}
\dfrac{\alpha_{i}^{1}(\mathbf{\mu}^{e}\mathbf{-}\frac{m_{2}}{M_{2}%
}\mathbf{\lambda}^{e}\mathbf{)}}{\mid\mathbf{\mu}^{e}\mathbf{-}\frac{m_{2}%
}{M_{2}}\mathbf{\lambda}^{e}\mid^{2i+3}}+\dfrac{Gm_{0}m_{1}}{M_{2}}%
{\displaystyle\sum\limits_{i=0}^{1}}
\dfrac{\alpha_{i}^{2}(\mathbf{\mu}^{e}\mathbf{+}\frac{m_{1}}{M_{2}%
}\mathbf{\lambda}^{e}\mathbf{)}}{\mid\mathbf{\mu}^{e}\mathbf{+}\frac{m_{1}%
}{M_{2}}\mathbf{\lambda}^{e}\mid^{2i+3}}\\
\\
(\mathbf{\nabla}_{\mathbf{\mu}}\mathcal{V})_{e}=Gm_{0}\left(
{\displaystyle\sum\limits_{i=0}^{1}}
\dfrac{\alpha_{i}^{1}(\mathbf{\mu}^{e}\mathbf{-}\frac{m_{2}}{M_{2}%
}\mathbf{\lambda}^{e}\mathbf{)}}{\mid\mathbf{\mu}^{e}\mathbf{-}\frac{m_{2}%
}{M_{2}}\mathbf{\lambda}^{e}\mid^{2i+3}}+%
{\displaystyle\sum\limits_{i=0}^{k}}
\dfrac{\alpha_{i}^{2}(\mathbf{\mu}^{e}\mathbf{+}\frac{m_{1}}{M_{2}%
}\mathbf{\lambda}^{e}\mathbf{)}}{\mid\mathbf{\mu}^{e}\mathbf{+}\frac{m_{1}%
}{M_{2}}\mathbf{\lambda}^{e}\mid^{2i+3}}\right)
\end{gather*}

where $\alpha_{0}^{1}=\alpha_{0}^{2}=m_{0}$, $\alpha_{1}^{1}=\beta
_{1}=3/2(C_{1}-A_{1})$ and $\alpha_{1}^{2}=\beta_{2}=3/2(C_{2}-A_{2}).$

The following identities are also verified%
\[
(\mathbf{\nabla}_{\mathbf{\lambda}}\mathcal{V})_{e}=(A_{11})_{e}%
\mathbf{\lambda}^{e}\mathbf{+(}A_{12})_{e}\mathbf{\mu}^{e}\mathbf{,\quad
}(\mathbf{\nabla}_{\mathbf{\mu}}\mathcal{V})_{e}=(A_{21})_{e}\mathbf{\lambda
}^{e}\mathbf{+(}A_{22})_{e}\mathbf{\mu}^{e}%
\]

where%
\begin{align*}
(A_{11})_{e}  &  =\dfrac{Gm_{1}m_{2}}{\mid\mathbf{\lambda}^{e}\mid^{3}}%
+\dfrac{Gm_{0}m_{2}^{2}}{M_{2}^{2}}\left(
{\displaystyle\sum\limits_{i=0}^{1}}
\dfrac{\alpha_{i}^{1}}{\mid\mathbf{\mu}^{e}\mathbf{-}\frac{m_{2}}{M_{2}%
}\mathbf{\lambda}^{e}\mid^{2i+3}}\right)  +\dfrac{Gm_{0}m_{1}^{2}}{M_{2}^{2}%
}\left(
{\displaystyle\sum\limits_{i=0}^{1}}
\dfrac{\alpha_{i}^{2}}{\mid\mathbf{\mu}^{e}\mathbf{+}\frac{m_{1}}{M_{2}%
}\mathbf{\lambda}^{e}\mid^{2i+3}}\right)  \bigskip\\
(A_{12})_{e}  &  =\dfrac{Gm_{0}m_{1}}{M_{2}}\left(
{\displaystyle\sum\limits_{i=0}^{1}}
\dfrac{\alpha_{i}^{2}}{\mid\mathbf{\mu}^{e}\mathbf{+}\frac{m_{1}}{M_{2}%
}\mathbf{\lambda}^{e}\mid^{2i+3}}\right)  -\dfrac{Gm_{0}m_{2}}{M_{2}}\left(
{\displaystyle\sum\limits_{i=0}^{1}}
\dfrac{\alpha_{i}^{1}}{\mid\mathbf{\mu}^{e}\mathbf{-}\frac{m_{2}}{M_{2}%
}\mathbf{\lambda}^{e}\mid^{2i+3}}\right)  \bigskip\\
(A_{12})_{e}  &  =(A_{21})_{e}\\
(A_{22})_{e}  &  =Gm_{0}\left(
{\displaystyle\sum\limits_{i=0}^{1}}
\dfrac{\alpha_{i}^{1}}{\mid\mathbf{\mu}^{e}\mathbf{-}\frac{m_{2}}{M_{2}%
}\mathbf{\lambda}^{e}\mid^{2i+3}}+%
{\displaystyle\sum\limits_{i=0}^{1}}
\dfrac{\alpha_{i}^{2}}{\mid\mathbf{\mu}^{e}\mathbf{+}\frac{m_{1}}{M_{2}%
}\mathbf{\lambda}^{e}\mid^{2i+3}}\right)
\end{align*}

\section{Some numerical results}

In order to obtain the values of $C_{i}-A_{i},$ $(i=1,2),$ the below
relationships will be utilized%
\begin{gather*}
C_{1}-A_{1}=(1-\mu)\left(  \dfrac{e_{S_{1}}}{Z}\right)  ^{2}J_{2}^{S_{1}}\\
\\
C_{2}-A_{2}=\mu\left(  \dfrac{e_{S_{2}}}{Z}\right)  ^{2}J_{2}^{S_{2}}%
\end{gather*}
where $e_{S_{i}}$ and $p_{S_{i}}$ represent the equatorial and polar radius of
$S_{i}$ $(i=1,2)$, $J_{2}^{S_{i}}=\frac{2}{5}\varepsilon_{i},$ respectively,
and $\varepsilon_{i}=\dfrac{e_{S_{i}}-p_{S_{i}}}{e_{S_{i}}}$. $S_{1},$ $S_{2}$
are considered to be homogeneous ellipsoids. The distances are measured in
kilometers.\bigskip

\begin{center}
\fbox{%
\begin{tabular}
[c]{llll}%
\underline{$S_{2}S_{1}S_{0}$}$\medskip(m_{0}\rightarrow0)$ & \underline
{Without oblat.} & \underline{Oblat. of $S_{2}$} & \underline{Oblat. of $S_{2}
$ and $S_{1}$}\\
Earth-Moon-$S_{0}\medskip\quad$ & $448879.206$ & $448879.221$ & $448879.251$\\
Mars-Phobos-$S_{0}\medskip$ & $9414.945$ & $9414.958$ & $9414.958$\\
\underline{$S_{0}S_{2}S_{1}$}$\medskip(m_{0}\rightarrow0)$ & \underline
{Without oblat.} & \underline{Oblat. of $S_{2}$} & \underline{Oblat. of $S_{2}
$ and $S_{1}$}\\
$S_{0}$-Earth-Moon$\medskip\quad$ & $381679.691$ & $381679.763$ & $381679.763
$\\
$S_{0}$-Mars-Phobos$\medskip$ & $9310.642$ & $9310.666$ & $9310.668$\\
\underline{$S_{2}S_{0}S_{1}$}$\medskip(m_{0}\rightarrow0)$ & \underline
{Without oblat.} & \underline{Oblat. of $S_{2}$} & \underline{Oblat. of $S_{2}
$ and $S_{1}$}\\
Earth-$S_{0}$-Moon$\medskip\quad$ & $326409.744$ & $326409.780$ & $326409.751
$\\
Mars-$S_{0}$-Phobos$\medskip$ & $9339.156$ & $9339.196$ & $9339.196$%
\end{tabular}
}
\end{center}

\section{Coefficients of the characteristic polynomial in Eulerian relative
equilibria $S_{0}S_{2}S_{1}$}

The coefficients of the characteristic polynomial (\ref{polinomiocarac}) are%
\begin{align*}
\omega_{e}^{2}  &  =\dfrac{G\,((m_{2}+m_{1})\,\rho^{4}+(2\,m_{1}%
+2\,m_{2})\,\rho^{3}+(m_{2}+m_{1})\,\rho^{2}-2\,m_{0}\,\rho-m_{0})}%
{{\lambda_{e}^{3}}(1+\rho)^{2}\rho^{2}}\bigskip\\
& \\
p  &  =\dfrac{G((m_{2}+4m_{0}+m_{1})\rho^{3}+(3m_{2}+6m_{0})\rho^{2}%
+(4m_{0}+3m_{2})\rho+m_{0}+m_{2})}{(1+\rho)^{3}\rho^{3}{\lambda_{e}^{3}}%
}\bigskip\\
& \\
q  &  =\dfrac{G((-2m_{1}\rho^{4}m_{2}+(-2m_{0}m_{1}+m_{1}^{2}+m_{2}^{2}%
-2m_{1}m_{2}-2m_{0}m_{2})\rho^{3}}{((1+\rho)^{3}\rho^{3}{\lambda_{e}^{3}})}+\\
&  \dfrac{(3m_{2}^{2}+m_{1}m_{2}-6m_{0}m_{1})\rho^{2}+(-m_{1}m_{2}+3m_{2}%
^{2}+2m_{0}m_{2}-4m_{0}m_{1})\rho}{((1+\rho)^{3}\rho^{3}{\lambda_{e}^{3}})}+\\
&  \dfrac{m_{2}^{2}-m_{0}m_{1}+m_{0}m_{2}-m_{1}m_{2}))}{((1+\rho)^{3}\rho
^{3}{\lambda_{e}^{3}})}\bigskip\\
& \\
r  &  =\dfrac{G^{2}(a_{1}\rho^{4}+a_{2}\rho^{4}+a_{3}\rho^{2}+a_{4}\rho
+a_{5})}{((1+\rho)^{8}\,\rho^{8}\,{\lambda_{e}^{9}})}%
\end{align*}

\subsection{Coefficients $a_{i}$ $(i=1,\ldots,5)$}

\begin{center}%
\begin{gather*}
a_{1}=-42m_{2}^{7}m_{1}-48m_{2}^{7}m_{0}-147m_{2}^{6}m_{1}^{2}-336m_{2}%
^{6}m_{1}m_{0}-129m_{2}^{6}m_{0}^{2}\smallskip\\
-207m_{2}^{5}m_{1}^{3}-782m_{2}^{5}m_{1}^{2}m_{0}-673m_{2}^{5}m_{1}m_{0}%
^{2}-81m_{2}^{5}m_{0}^{3}-150m_{2}^{4}m_{1}^{4}\smallskip\\
-869m_{2}^{4}m_{1}^{3}m_{0}-1325m_{2}^{4}m_{1}^{2}m_{0}^{2}-378m_{2}^{4}%
m_{1}m_{0}^{4}-64m_{2}^{3}m_{1}^{5}\smallskip\\
-513m_{2}^{3}m_{1}^{4}m_{0}-1270m_{2}^{3}m_{1}^{3}m_{0}^{2}-702m_{2}^{3}%
m_{1}^{2}m_{0}^{3}-14m_{2}^{2}m_{1}^{6}\smallskip\\
-165m_{2}^{2}m_{1}^{5}m_{0}-610m_{2}^{2}m_{1}^{4}m_{0}^{2}-648m_{2}^{2}%
m_{1}^{3}m_{0}^{3}-24m_{2}m_{1}^{6}m_{0}\smallskip\\
-119m_{2}m_{1}^{5}m_{0}^{2}-297m_{2}m_{1}^{4}m_{0}^{3}+2m_{1}^{6}m_{0}%
^{2}-54m_{1}^{5}m_{0}^{3}\smallskip
\end{gather*}%
\begin{gather*}
a_{2}=-60m_{2}^{7}m_{1}-54m_{2}^{7}m_{0}-243m_{2}^{6}m_{1}^{2}-474m_{2}%
^{6}m_{1}m_{0}-173m_{2}^{6}m_{0}^{2}\smallskip\\
-399m_{2}^{5}m_{1}^{3}-1345m_{2}^{5}m_{1}^{2}m_{0}-999m_{2}^{5}m_{1}m_{0}%
^{2}-135m_{2}^{5}m_{0}^{3}-329m_{2}^{4}m_{1}^{4}\smallskip\\
-1846m_{2}^{4}m_{1}^{3}m_{0}-2223m_{2}^{4}m_{1}^{2}m_{0}^{2}-648m_{2}^{4}%
m_{1}m_{0\ }^{3}-138m_{2}^{3}m_{1}^{5}\smallskip\\
-1364m_{2}^{3}m_{1}^{4}m_{0}-2506m_{2}^{3}m_{1}^{3}m_{0}^{2}-1242m_{2}%
^{3}m_{1}^{2}m_{0}^{3}-24m_{2}^{2}m_{1}^{6}\smallskip\\
-536m_{2}^{2}m_{1}^{5}m_{0}-1530m_{2}^{2}m_{1}^{4}m_{0}^{2}-1188m_{2}^{2}%
m_{1}^{3}m_{0}^{3}-90m_{2}m_{1}^{6}m_{0}\smallskip\\
-477\,m_{2}m_{1}^{5}m_{0}^{2}-567m_{2}m_{1}^{4}m_{0}^{3}-56m_{1}^{6}m_{0}%
^{2}-108m_{1}^{5}m_{0}^{3}\smallskip
\end{gather*}%
\begin{gather*}
a_{3}=-42m_{2}^{7}m_{1}-36m_{2}^{7}m_{0}-183m_{2}^{6}m_{1}^{2}-342m_{2}%
^{6}m_{1}m_{0}-93m_{2}^{6}m_{0}^{2}\smallskip\\
-349m_{2}^{5}m_{1}^{3}-1097m_{2}^{5}m_{1}^{2}m_{0}-630m_{2}^{5}m_{1}m_{0}%
^{2}-81m_{2}^{5}m_{0}^{3}-358m_{2}^{4}m_{1}^{4}\smallskip\\
-1776m_{2}^{4}m_{1}^{3}m_{0}-166\,m_{2}^{4}m_{1}^{2}m_{0}^{2}-405m_{2}%
^{4}m_{1}m_{0}^{3}-189m_{2}^{3}m_{1}^{5}\smallskip\\
-1614m_{2}^{3}m_{1}^{4}m_{0}-2256m_{2}^{3}m_{1}^{3}m_{0}^{2}-810m_{2}^{3}%
m_{1}^{2}m_{0}^{3}-31m_{2}^{2}m_{1}^{6}\smallskip\\
-827m_{2}^{2}m_{1}^{5}m_{0}-1683m_{2}^{2}m_{1}^{4}m_{0}^{2}-810m_{2}^{2}%
m_{1}^{3}m_{0}^{3}-6m_{2}\,m_{1}^{7}\smallskip\\
-228m_{2}m_{1}^{6}m_{0}-666m_{2}m_{1}^{5}m_{0}^{2}-405m_{2}m_{1}^{4}m_{0}%
^{3}-30m_{1}^{7}m_{0}\smallskip\\
-81m_{1}^{5}m_{0}^{3}-111m_{1}^{6}m_{0}^{2}\smallskip
\end{gather*}%
\begin{gather*}
a_{4}=-12m_{2}^{7}m_{1}-12m_{2}^{7}m_{0}-56m_{2}^{6}m_{1}^{2}-114m_{2}%
^{6}m_{1}m_{0}-24m_{2}^{6}m_{0}^{2}\smallskip\\
-130m_{2}^{5}m_{1}^{3}-387m_{2}^{5}m_{1}^{2}m_{0}-162m_{2}^{5}m_{1}m_{0}%
^{2}-179m_{2}^{4}m_{1}^{4}\smallskip\\
-687m_{2}^{4}m_{1}^{3}m_{0}-432m_{2}^{4}m_{1}^{2}m_{0}^{2}-140m_{2}^{3}%
m_{1}^{5}-693m_{2}^{3}m_{1}^{4}m_{0}\smallskip\\
-588m_{2}^{3}m_{1}^{3}m_{0}^{2}-52m_{2}^{2}m_{1}^{6}-387m_{2}^{2}m_{1}%
^{5}m_{0}-432m_{2}^{2}m_{1}^{4}m_{0}^{2}-6m_{2}m_{1}^{7}\smallskip\\
-108m_{2}m_{1}^{6}m_{0}-162m_{2}m_{1}^{5}m_{0}^{2}-12m_{1}^{7}m_{0}%
-24m_{1}^{6}m_{0}^{2}\smallskip
\end{gather*}%
\begin{gather*}
a_{5}=-(m_{0}+m_{2})(18m_{0}m_{2}^{6}+12m_{1}m_{2}^{6}+94m_{2}^{5}m_{0}%
m_{1}+36m_{1}^{2}m_{2}^{5}\smallskip\\
+81m_{2}^{4}m_{0}^{2}m_{1}+168m_{2}^{4}m_{0}m_{1}^{2}+42m_{2}^{4}m_{1}%
^{3}+128m_{2}^{3}m_{0}m_{1}^{3}\smallskip\\
+27m_{2}^{3}m_{1}^{4}+15m_{2}^{2}m_{1}^{5}+31m_{2}^{2}m_{0}m_{1}^{4}%
+126m_{2}^{2}m_{0}^{2}m_{1}^{3}+18m_{0}^{2}m_{2}^{5}\smallskip\\
+54m_{2}m_{0}^{2}m_{1}^{4}+12m_{2}m_{0}m_{1}^{5}+5m_{2}m_{1}^{6}+7m_{1}%
^{6}m_{0}+9m_{0}^{2}m_{1}^{5}\smallskip\\
+144m_{2}^{3}m_{0}^{2}m_{1}^{2})\smallskip
\end{gather*}

\textbf{Acknowledgements}
\end{center}

This research was partially supported by the Spanish Ministerio de Ciencia y
Tecnolog\'{\i}a (Project BFM2003-02137) and by the Consejer\'{\i}a de
Educaci\'{o}n y Cultura de la Comunidad Aut\'{o}noma de la Regi\'{o}n de
Murcia (Project S\'{e}neca 2002: PC-MC/3/00074/FS/02).

\end{document}